\documentclass{elsarticle}

\usepackage{amssymb,amsmath}
\usepackage{graphicx,overpic,psfrag}
\usepackage{a4wide}
\usepackage{subfigure}
\usepackage[arrow, matrix, curve,cmtip]{xy}
\usepackage{pb-diagram}
\usepackage{color}

\newtheorem{lemma}{Lemma}
\newtheorem{theorem}[lemma]{Theorem}
\newtheorem{corollary}[lemma]{Corollary}
\newtheorem{definition}[lemma]{Definition}

\newtheorem{proposition}[lemma]{Proposition}
\newtheorem{remark}{Remark}

\def\be{\begin{eqnarray}}
\def\ee{\end{eqnarray}}
\def\bes{\begin{eqnarray*}}
\def\ees{\end{eqnarray*}}

\def\R{\mathbb R}       
\def\P{\mathbb P}       

\def\av{\mbox{\bf a}} \def\bv{\mbox{\bf b}} \def\cv{\mbox{\bf c}}
\def\dv{\mbox{\bf d}} \def\ev{\mbox{\bf e}} \def\fv{\mbox{\bf f}}
\def\gv{\mbox{\bf g}}  
  
\def\mv{\mbox{\bf m}} \def\nv{\mbox{\bf n}} 
\def\pv{\mbox{\bf p}} \def\qv{\mbox{\bf q}} 
\def\sv{\mbox{\bf s}}  \def\uv{\mbox{\bf u}}
  \def\xv{\mbox{\bf x}}
\def\yv{\mbox{\bf y}} \def\zv{\mbox{\bf z}}

\def\rk{\operatorname{rk~}}     
\def\det{\operatorname{det}}    
\def\diag{\operatorname{diag}}  

\def\eps{\varepsilon}

\def\ol{\overline}
\def\wt{\widetilde}

\def\phi{\varphi}

\begin{document}
\begin{frontmatter}

\title{The Relation Between Offset and Conchoid Constructions}


\author[mp]{Martin Peternell}
\author[mp]{Lukas Gotthart}
\author[rs]{J.~Rafael Sendra}
\author[js]{Juana Sendra}

\address[mp]{Institute of Discrete Mathematics and Geometry, 
Vienna University of Technology, Vienna, Austria}
\address[rs]{Dpto.~de Fisica y Matem\'aticas,
Universidad de Alcal\'a, E-28871 Alcal\'a de Henares, Madrid, Spain}
\address[js]{Dpto.~Matem\'atica Aplicada a la I.T.~de Telecomunicaci\'on. 
Research Center on Software Technologies and Multimedia Systems 
for Sustainability (CITSEM), UPM, Spain}

\begin{abstract}

The one-sided offset surface $F_d$ of a given surface $F$ is, roughly speaking,
obtained by shifting the tangent planes of $F$ in direction of its oriented
normal vector. The conchoid surface $G_d$ of a given surface $G$ is roughly speaking
obtained by increasing the distance of $G$ to a fixed reference point $O$ by $d$.
Whereas the offset operation is well known and implemented in most CAD-software systems,
the conchoid operation is less known, although already mentioned by the ancient Greeks,
and recently studied by some authors.

These two operations are algebraic and create new objects from given input objects.
There is a surprisingly simple relation between the offset and the conchoid operation.
As derived there exists a rational bijective quadratic map which transforms a given surface
$F$ and its offset surfaces $F_d$ to a surface $G$ and its conchoidal surfaces $G_d$,
and vice versa. Geometric properties of this map are studied and illustrated at hand 
of some examples. Furthermore, rational universal parameterizations for
offsets and conchoid surfaces are provided.
\end{abstract}


\end{frontmatter}

\section{Introduction}\label{intro:sec}

There is a large variety of contributions dealing with the geometry of offsets constructions discussing different aspects. Since we are focusing on parametrization problems 
we mention here only some of them, see for instance \cite{arrondo}, \cite{F1}, 
\cite{F2}, \cite{KPet}, \cite{pet1}, \cite{pot98}, \cite{sendra0}, \cite{VL} and 
references on the topic in~\cite{libro}. Conchoidal constructions, 
although not so extensively studied, 
have been recently addressed by different authors too, see for 
instance~\cite{alberto}, \cite{pet10}, \cite{pet12}, \cite{sendra1}, \cite{sendra2},
\cite{sendra3}, \cite{sendra4}. 
Both geometric constructions were already utilized in the past 
(Leibnitz studied parallel curves and ancient Greeks used conchoids), 
and nowadays are used in practical applications 
(see e.g. \cite{F1}, \cite{F2}, \cite{KPet}, \cite{HL} and introduction of \cite{sendra1}).

\paragraph{Contribution:}
Considering an algebraic irreducible surface $F$ and its continuous family of 
one-sided offset surfaces $F_d$, for $d\in\R$,
there exists a birational quadratic map $\alpha$ so that the surface $G_d=\alpha(F_d)$
is the one-sided conchoid surface of $G=\alpha(F)$ for $d\in\R$. The inverse
map realizes the correspondence between a family of one-sided conchoid surfaces $G_d$
and a family of one-sided offset surfaces $F_d$.
Since $\alpha$ is a birational map, rationality is preserved in both directions.
All geometric properties and results, which are known for one family can be
transformed to properties and results for the other family.
To derive this correspondence, one-sided offset surfaces are considered as envelopes of
tangent planes. In addition, we introduce $\alpha$ for implicitly defined algebraic surfaces, 
and thus the results extend to the two-sided offsets and conchoids.
Throughout this paper, we present the results for surfaces but they are valid 
for hypersurfaces, in particular for plane algebraic curves.

Sections~\ref{intro:sec} and~\ref{offset:sec} introduce to representations
of offset surfaces and conchoid surfaces. Section~\ref{footpnt:sec} discusses the
foot-point map realizing the correspondence between offsets and conchoids. Further,
Sections~\ref{elemexam:sec} and~\ref{addexam:sec} present elementary examples
and with offsets and conchoids of ruled surfaces and quadrics more advanced ones.
Finally the conclusion also outlines the relations to bisector surfaces.
All computations and figures are carried out with the mathematical software Maple.

\paragraph{Notation:}

The scalar product of two vectors $\xv,\yv\in \R^3$  
is denoted as $\xv\cdot\yv$. Points in the Euclidean space $\R^3$ are identified 
with their coordinate vectors $\xv=(x,y,z)\in\R^3$ with respect to a chosen coordinate system. 
The projective extension of the Euclidean space $\R^3$ is denoted as $\P^3$.
Points $X$ in projective space $\P^3$ are identified with their homogeneous 
coordinate vectors $X=(x_0,x_1,x_2,x_3)\R$ $=$ $(x_0:x_1:x_2:x_3)$.
Let $\omega=\P^3\setminus\R^3$ be the plane at infinity in $\P^3$, 
that is the plane defined by $x_0=0$. 
Ideal points, that are points at infinity, are represented by $(0,\xv)\R$, with $\xv\in\R^3$.
For points $X\not\in\omega$, thus $x_0\not=0$, the relation between
Cartesian and homogeneous coordinates is given by
\bes
x = \frac{x_1}{x_0}, y = \frac{x_2}{x_0}, z = \frac{x_3}{x_0}.
\ees
Planes $U$ in $\P^3$ are the zero-set of a linear equation $U: u_0x_0 +\ldots+ u_3x_3=0$.
The coefficients $(u_0, \ldots, u_3)$ are the homogeneous coordinates of $U$. 
Identifying $U$ with these coordinates, we write $U=\R(u_0,u_1,u_2,u_3)$.
The ideal plane $\omega$ is represented by $\R(1,0,0,0)$.
To distinguish points from planes, we use $(x_0, \ldots, x_3)\R$ for points 
and $\R(u_0,\ldots,u_3)$ for planes. In this way a plane $U\subset \P^3$
is identified with a point $U=\R(u_0,u_1,u_2,u_3)\in{\P^3}^\star$, where
${\P^3}^\star$ denotes the dual space associated to $\P^3$.



\subsection{Offset construction}\label{offset_intro:sec}

Consider an irreducible algebraic surface $F\subset\R^3$ given by a (non necessarily rational) parametrization $\fv(u,v)$. Let $\nv(u,v)$ be the unit normal vector field of $\fv(u,v)$.
Then the {\sf one-sided offset surface} $F_d$ of $F$ at oriented distance $d\in \R$ 
is given by the parametrization
\begin{equation}\label{offset_surf:eq}
\fv_d(u,v) = \fv(u,v) + d\nv(u,v), \mbox{ with } \|\nv \|=1.
\end{equation}
Since $\nv$ is normalized, its partial derivatives $\nv_u$ and $\nv_v$
are orthogonal to $\nv$. This implies that the tangent planes of the offset
surface $F_d$ are parallel to the tangent planes of $F$ at corresponding
points $\fv(u,v)$ and $\fv_d(u,v)$.

Often the offset surface of the surface $F$ is defined as the envelope of a family 
of spheres of radius $d$, centered at the base surface $F$.
Let us denote this offset as $O_d(F)$. This definition obviously differs
from the definition given by~\eqref{offset_surf:eq}. 
Let us analyze the relation between these two concepts. 
In~\cite{sendra0} it is shown that $O_d(F)$ is algebraic and that it has 
at most two irreducible components. 
If $O_d(F)$ is reducible then the two components are 
$F_d$ and $F_{-d}$, the one-sided offsets at distances $d$ and $-d$. 
If $O_d(F)$ is irreducible, then $O_d(F)$ is the Zariski closure of $F_d$ and 
also of $F_{-d}$. In any case, $O_d(F)=F_d \cup F_{-d}$. 
But it should be noted that, if $O_d(F)$ is irreducible, then $F_d $ and $F_{-d}$ 
are not algebraic surfaces in the strict sense, since they 
represent the exterior and interior analytic components of the two-sided offset $O_d(F)$. 
Focusing on rational surfaces, we have the following definition.

\begin{definition}\label{ratoffset:def}
A rational surface $F$ is called {\sf rational offset surface} if $F$ admits 
a rational parametrization $\fv(u,v)$ with rational unit normal vector field $\nv(u,v)$ 
of $\fv(u,v)$.
\end{definition}

Note that if $F$ is a rational offset surface, then both $F_d$ and $F_{-d}$ admit 
a rational representation~\eqref{offset_surf:eq}. Furthermore,
\begin{itemize}
\item if $O_d(F)$ is irreducible then $O_d(F)$  is rational iff $F_d$ admits a rational representation~\eqref{offset_surf:eq}, and analogously for $F_{-d}$; 
\item if $O_d(F)$ is reducible then all components of $O_d(F)$  are rational iff $F_d$ and 
$F_{-d}$ admit a rational representation~\eqref{offset_surf:eq}.
\end{itemize}


\subsection{Conchoidal construction}\label{conchoid_constr:sec}

Consider, as above, an irreducible algebraic surface $G\subset \R^3$, $d\in \R$, and a fixed reference point $O$. Without loss of generality we assume that $O$ is the origin of 
a Cartesian coordinate system. Let $G$ be represented by a polar 
(not necessarily rational) representation
$\gv(u,v)=r(u,v)\sv(u,v)$, with $\|\sv(u,v)\|=1$. We call $\sv(u,v)$ the
{\sf spherical part} of the parameterization $\gv(u,v)$ and $r(u,v)$ the {\sf radius function}. 
In this situation,
the {\sf one-sided conchoid surface} $G_d$ of $G$ is obtained by increasing the radius
$r(u,v)$ by $d$ and thus $G_d$ admits the polar representation
\begin{equation}\label{conchsurf_par:eqn}
\gv_d(u,v)=(r(u,v) + d)\sv(u,v).
\end{equation}

In \cite{alberto}, \cite{sendra1}, \cite{sendra2}, the conchoidal variety is 
introduced considering both analytic components $G_d$ and $G_{-d}$. 
More precisely, the conchoid surface $C_d(G)$ of $G$ with respect to $O$ at distance $d$ 
is the (Zariski closure) set of points $Q$ in the line $OP$ 
at distance $d$ of a moving point $P\in G$. That is the Zariski closure of the set
\begin{equation}\label{conchsurf_geom:eqn}
C_d(G)=\{ Q\in OP \mbox{ with } P\in G, \mbox{ and } \overline{QP}=d\}.
\end{equation}

As in the case of offset surfaces in Section~\ref{offset_intro:sec} the two given 
definitions~\eqref{conchsurf_par:eqn} and~\eqref{conchsurf_geom:eqn} 
for conchoids differ, but are clearly related.
The conchoid surface $C_d(G)$ has at most two irreducible components (see~\cite{sendra1}, \cite{sendra2}). If $C_d(G)$ is reducible then the two components are $C_d(G)$ and $C_{-d}(G)$, 
the one-sided conchoids at distances $d$ and $-d$. If $C_d(G)$ is irreducible, then $C_d(G)$ is  the Zariski closure of $G_d$, and also of $G_{-d}$. 
In any case, $C_d(F)=G_d \cup G_{-d}$. But it should be noted that, if $C_d(G)$ is 
irreducible, then $G_d $ and $G_{-d}$ are not algebraic surfaces in the strict sense, 
since they represent the exterior and interior analytic components of 
the two-sided conchoid $C_d(F)$. 
To deal with rational surfaces $G$ and their conchoid surfaces $G_d$,
we define the following.
\begin{definition}\label{ratconch:def}
A surface $G$ is called {\sf rational conchoid surface} with respect to
the reference point $O$ if $G$ admits a {\sf rational polar representation} $r(u,v)\sv(u,v)$, 
with a rational radius function $r(u,v)$ and a rational parametrization $\sv(u,v)$ of 
the unit sphere $S^2$.
\end{definition}
Note that if $G$ is a rational conchoid surface, $G_d$ admits 
the rational representation~\eqref{conchsurf_par:eqn}. 
Furthermore, an analogous remark as on the rationality of the offset done above, 
is valid also for conchoids. 

\section{Representation of Offset and Conchoid surfaces}\label{offset:sec}

This section introduces the representations of the offset and conchoid surfaces.
While conchoid surfaces are represented as point sets, offset surfaces are represented
as envelopes of tangent planes. Performing in this way, we realize that their 
representations~\eqref{offset:rep} and~\eqref{conch:rep} are closely related. 

\subsection{Dual representation of offset surfaces}\label{dualrep_offset:sec}

Consider an irreducible algebraic surface $F\subset\R^3$, and its offset surfaces
$F_d\subset\R^3$. Let $\fv(u,v)$ be an affine parameterization of $F$, and
its offsets are thus represented by~\eqref{offset_surf:eq}. Dealing with offset
surfaces in this context, it is preferable to consider $F$ as envelope of its
tangent planes 
\begin{equation}\label{offset_surf_dual:eq}
E(u,v): \nv(u,v)\cdot\xv = e(u,v), \mbox{ with } e=\fv\cdot\nv,
\end{equation}
where $\nv(u,v)$ denotes a unit normal vector field of $\fv(u,v)$.
The function $e(u,v)$ represents the oriented distance of the
origin $O$ from the planes $E(u,v)$. 
Consequently the offset surfaces $F_d$ of $F$ have tangent planes
\begin{equation}\label{offset_surf_dual2:eq}
E_d(u,v): \nv(u,v)\cdot\xv = e(u,v)+d.
\end{equation}
Let $\ol{F}\subset \P^3$ be the projective surface corresponding to $F$. 
The tangent planes $E$ of $\ol F$ determine a surface $\ol F^\star\subset {\P^3}^\star$,
called the {\sf dual surface} of $\ol F$. 
While a parameterization of $\ol F$ in terms
of homogeneous point coordinates is $(1,\fv(u,v))\R$, a parameterization of 
$\ol F^\star$ and its offset surfaces $\ol F_d^\star$ reads
\bes
E = \R(-e(u,v), \nv(u,v)), \mbox{ and } E_d = \R(-e(u,v)-d, \nv(u,v)),
\ees
where we use same symbols for the planes $E$, $E_d$ and their homogeneous
coordinate vectors. The homogeneous coordinates of planes are the coefficients
of the defining linear polynomial. 

Given a two parameter family of planes $E(u,v): \nv(u,v)\cdot\xv = e(u,v)$, 
an affine parameterization $\fv(u,v)$ of the envelope surface $F$ is computed 
as solution of the system of linear equations
\be\label{envelope_cond:eq}
\{ \nv(u,v)\cdot\xv = e(u,v),  \,
\nv_u(u,v)\cdot\xv = e_{u}(u,v), \,
\nv_v(u,v)\cdot\xv = e_{v}(u,v)
\}
\ee
where $X_u, X_v$ denote the partial derivatives of $X(u,v)$
with respect to $u$ and $v$.
Concerning the envelope surface $F$ of planes $E(u,v)$ there exist some
degenerate cases, depending on the rank of the coefficient matrix $(\nv, \nv_u, \nv_v)$.
If $\nv$ and $e$ are constant, $F$ is a plane. If $\nv$ and $e$ are functions of one
variable only, $\nv$ represents a curve in the unit sphere $S^2$, and 
consequently the surface $F$ is developable. 
Otherwise if $\nv(u,v)$ represents a surface in $S^2$ and $F$ is typically a 
non-developable surface in $\R^3$. But even if $\nv(u,v)$ is two-dimensional, 
it might happen that the envelope $F$ of the planes $E$ is a single point, 
see Section~\ref{sphere_offsets}, or that the planes $E$ are tangent planes 
of a curve, see Section~\ref{pedal_conics:sec}.

Let $\ol F\subset\P^3$ be an irreducible algebraic surface, defined as zero-set of a homogeneous
polynomial $\ol F(x_0,x_1,x_2,x_3)$. To obtain a more compact notation, we use the 
notation $\ol{F}(x_0,\xv)$, with $\xv=(x_1,x_2,x_3)$. 
Likewise, the dual surface $\ol{F}^\star\subset {\P^3}^\star$ is the zero-set
of a homogeneous polynomial $\ol F^\star(u_0, \uv)$, with $\uv=(u_1,u_2,u_3)$.
Since the tangent planes of $\ol F(x_0,\xv)=0$ are represented by partial derivatives
$\ol F_{x_i}$ of $\ol F$ with respect to $x_i$, 
the homogeneous plane coordinates of $\ol F$ are
\be\label{gradient_map}
(u_0: u_1: u_2: u_3) = (\ol F_{x_0}, \ldots, \ol F_{x_3})(x_0, \ldots, x_3) = 
\nabla(\ol{F})(x_0,\xv), 
\ee
evaluated at points $X=(x_0:\ldots:x_3)\in \ol F$.   
The implicit equation $\ol{F}^\star(u_0,\uv)$ of the dual surface $\ol F^\star$ 
can be computed by eliminating $x_0,\xv$ in 
the algebraic system $\{\ol{F}(x_0,\xv)=0,(u_0,\uv)=\nabla(\ol{F})(x_0,\xv)\}$.

Let us assume that $\ol{F}$ is of class $n$, which means that 
the dual surface $\ol{F}^\star$ has degree $n$. 
The class of $\ol{F}$ expresses the algebraically counted number of 
tangent planes passing through a generic line in $\P^3$, whereas the degree of 
a surface $\ol{F}$ counts the intersection points with a generic line.
Let $\ol{F}^{\star}$ be expressed as
\be\label{fstar_implicit}
\ol{F}^\star(u_0, \uv) = u_0^nf_0 + u_0^{n-1}f_1(\uv) + \ldots + u_0^{n-j}f_j(\uv)+
\ldots + u_0f_{n-1}(\uv) + f_n(\uv),
\ee
where  $f_j(\uv)$ are homogeneous polynomials of degree $j$ in $\uv$. 
The plane at infinity $\omega: x_0=0$
is tangent to $\ol{F}$, or equivalently $\omega\in \ol{F}^\star$,
if $f_0=0$. Furthermore, $\omega$ is an $r$-fold plane of $\ol{F}^\star$, 
exactly if $f_0=\cdots =f_{r-1}=0$, but $f_r\not=0$.

\subsection{Parametric representation of offset and conchoid surfaces}
\label{offset_conchoid_parameterization:sec}

Let $\cal P$ be the point set of $\R^3$ and let $\cal E$ be the set of planes 
in $\R^3$. Considering the offset construction, an irreducible affine surface $F$ 
is defined as envelope of its tangent planes~\eqref{offset_surf_dual:eq}. 
Let $S^2: x^2+y^2+z^2=1$ be the unit sphere in $\R^3$, then the tangent planes $E$ 
of $F$ are defined by the map
\begin{eqnarray}\label{offset:rep}
\phi: S^2\times\R  &\to& \cal E \nonumber\\
      (\nv, e) & \mapsto & E:\nv\cdot\xv=e,
\mbox{ with } \nv\in S^2, e\in \R.
\end{eqnarray}
Any rational parametrization $\nv(u,v)$ of $S^2$ and any rational radius function $e(u,v)$
define a rational offset surface $F$ in the sense of Definition~\ref{ratoffset:def}.
Since the offset surface $F_d$ is obtained as envelope of 
the planes~\eqref{offset_surf_dual2:eq}, the offset map $o^\star$ with respect to
the offset distance $d$ is defined by
\be\label{offset_map:eq}
o^\star: {\cal E} &\to& {\cal E}  \nonumber\\
         \phi(\nv, e) &\mapsto& \phi(\nv, e+d).
\ee 
Considering the conchoid construction,
an irreducible algebraic surface $G\subset \R^3$ is represented by a polar representation $\gv(u,v)=r(u,v)\sv(u,v)$,
with $\|\sv(u,v) \|=1$. Thus $G$ is defined by the map
\begin{eqnarray}\label{conch:rep}
\gamma: S^2\times\R  &\to& \cal P \nonumber\\
        (\sv, r) & \mapsto & \gv = r\sv,
        \mbox{ with } \sv\in S^2, r\in \R.
\end{eqnarray}
Any rational parametrization $\sv(u,v)$ of $S^2$ and any rational radius function $r(u,v)$
define a rational conchoid surface $G$ in the sense of Definition~\ref{ratconch:def}.
The conchoidal map $c$, specifying the relation
between a surface $G$ and its conchoid surface $G_d$, is defined as
\be\label{conchoid_map:eq}
c: {\cal P} &\to& {\cal P}, \nonumber\\
    \gamma(\sv, r) &\mapsto& \gamma(\sv, r+d).
\ee 

We note that $\phi$ as well as $\gamma$ are considered as local parameterizations
of the tangent planes or the points of surfaces in $\R^3$. It is obvious that 
$\phi(\nv,e)$ and $\phi(-\nv,-e)$ define the same non-oriented plane in $\cal E$.
Likewise $\gamma(\sv,r)$ and $\gamma(-\sv,-r)$ define the same point in $\cal P$.
Thus, $\phi$ and $\gamma$ are not injective maps, but any element in $\cal E$ and
$\cal P$, respectively, has two pre-images. 
If injectivity of the maps is an issue, one can identify antipodal points in $S^2\times \R$.
Another possibility to overcome this problem for the map $\phi$ is by defining 
oriented planes in $\R^3$, 
where the orientation is determined by the oriented normal vector $\nv$.
When having introduced the one-sided offset $F_d$ in~\eqref{offset_surf:eq}, 
we have used an oriented normal vector field. But when dealing with algebraic objects,
the offset surface $O_d(F)$ contains both, the inner and outer offset. The same holds
for the conchoid surface $C_d(G)$ of an algebraic surface $G$.
Since the construction in Section~\ref{footpnt:sec} ignores orientations, 
we consider planes in $\R^3$ as non-oriented. Additionally we notice that when 
parameterizing surfaces and constructing offsets or conchoids, the base surface 
is traced twice, in order to represent both components of the offset or conchoid.

\subsection{Rational parameterizations of $S^2$}
\label{ratpara_s2:sec}

Since the focus is on rational families of offset surfaces
and conchoidal surfaces we discuss universal rational parameterizations of $S^2$.
Following \cite{dietz} we choose four arbitrary rational functions
$a(u,v)$, $b(u,v)$, $c(u,v)$ and $d(u,v)$ without common factor.
Let
\bes
A = 2(ac+bd), B = 2(bc-ad), C = a^2+b^2-c^2-d^2,
D = a^2+b^2+c^2+d^2,
\ees
then $\qv(u,v) = \frac{1}{D}(A,B,C)$ is a rational parametrization of 
the unit sphere $S^2$. Thus $\phi(\qv(u,v), \rho(u,v))$ with a rational 
function $\rho(u,v)$ defines a rational parametrization of a rational 
offset surface and likewise $\gamma(\qv(u,v), \rho(u,v))$ defines 
a rational parametrization of a rational conchoid surface. 
The similarities between representations~\eqref{offset:rep} and~\eqref{conch:rep}
indicate that there is a close relation between the offset construction
and the conchoid construction. 
These correspondences are studied in Section~\ref{footpnt:sec}.


\section{The foot-point map}\label{footpnt:sec}

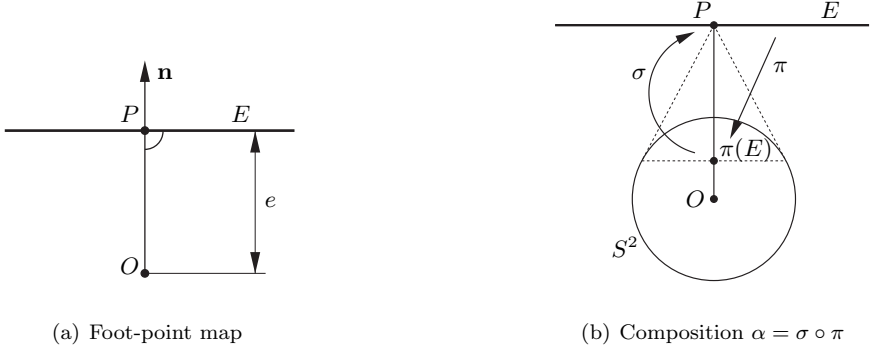
\begin{figure}[ht]
\subfigure[Foot-point map]{\label{footpoint_fig}
\begin{overpic}[width=.3\textwidth]{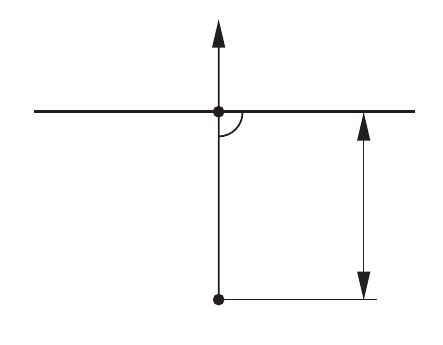}%
{\small
\put(75,55){$E$} \put(42,55){$P$} \put(42,10){$O$}
\put(85,30){$e$}  \put(53,68){$\nv$}
}
\end{overpic}
}
\hfil
\subfigure[Composition $\alpha=\sigma\circ\pi$]{\label{decompose_alpha_fig}
\begin{overpic}[width=.3\textwidth]{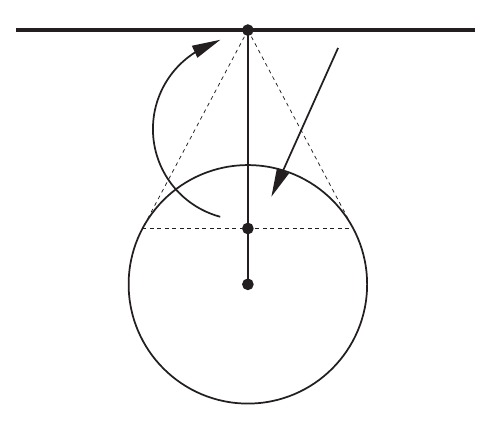}%
{\small
\put(42,30){$O$}
\put(82,86){$E$} \put(44,86){$P$} \put(52,45){$\pi(E)$}
\put(20,15){$S^2$} \put(68,70){$\pi$} \put(26,68){$\sigma$}
}
\end{overpic}
}
\caption{Geometric properties of the foot-point map $\alpha$}
\label{footpoint_map:fig}
\end{figure}

According to Section~\ref{conchoid_constr:sec}, we consider the origin $O$ 
in $\R^3$ as reference point for the conchoidal construction. We now introduce a map, 
that establishes the connection between offsets and conchoids.
The foot-point map $\alpha$ with respect to $O$ is defined as
\begin{eqnarray} \label{footpoint_map:eq}
\alpha: {\cal E} &\to& {\cal P}\\
    E: \xv\cdot\nv=e &\mapsto& P=\alpha(E)=\frac{e}{\|\nv\|^2}\nv, \nonumber
\end{eqnarray}
and maps planes $E\in \cal E$ to points $P\in \cal P$ of $\R^3$, see
Figure~\ref{footpoint_fig}.
The map $\alpha$ is rational and bijective except for planes $E$ passing through $O$.
The inverse map $\alpha^{-1}$ equals the dual map $\alpha^\star$, and 
transforms points $P\not=O$ to planes $E=\alpha^{-1}(P)$ which have $OP$ as normal and 
pass through $P$.
It reads
\begin{eqnarray} \label{footpoint_map_inv:eq}
\alpha^\star: {\cal P}\setminus \{ O \} &\to& {\cal E}\\
    P= \pv &\mapsto& \alpha^\star(P) = E: \xv\cdot\pv=\pv\cdot\pv, \nonumber
\end{eqnarray}
The maps $\alpha$ and $\alpha^\star$ are the basic ingredients
to construct rational conchoid surfaces from rational offset surfaces and vice versa. 
If the reference point $O$ is replaced by some other point $Z$, the foot-point map 
$\alpha_Z$ with respect to $Z$ reads
\begin{eqnarray*}
\alpha_Z: {\cal E} &\to& {\cal P}\\
E: \xv^T\cdot\nv=e &\mapsto& P=\alpha_Z(E)=\zv+\frac{e-\zv\cdot \nv}{\|\nv\|^2}\,\nv, 
\nonumber
\end{eqnarray*}

The map~\eqref{footpoint_map:eq} is a quadratic
plane-to-point mapping. To get more insight to the correspondence between
offsets and conchoids we provide representations of $\alpha$ and $\alpha^\star$
in terms of homogeneous coordinates.
Consider the projective extension $\P^3$ of $\R^3$, and its dual space ${\P^3}^\star$. 
The points $\R(u_0,\ldots,u_3)\in{\P^3}^\star$ are identified with the planes 
$U\subset P^3$.

The map $\alpha:{\cal E}\to{\cal P}$ uniquely determines the map 
$\ol{\alpha}: {\P^3}^\star \to \P^3$, and for simplicity both maps 
are denoted by $\alpha$. In terms of homogeneous 
coordinates it reads  
\be \label{footpoint_map2:eq} 
\begin{array}{rlll}
\alpha: {\P^3}^\star & \to    & \P^3&  \\
U=\R(u_0,\ldots,u_3) &\mapsto & X    &= (x_0,x_1,x_2,x_3)\R \\
& & &=(-(u_1^2+u_2^2+u_3^2),u_0u_1,u_0u_2,u_0u_3)\R.\\
\end{array}
\ee
The exceptional planes of $\alpha$ are the ideal plane $\omega: x_0=0$, and the 
tangent planes of the {\sf isotropic cone} $\Delta: u_0=0, u_1^2+u_2^2+u_3^2=0$ 
with vertex at $O=(1,0,0,0)\R$.
The dual map $\alpha^\star=\alpha^{-1}$ reads 
\be \label{footpoint_map3:eq}
\begin{array}{rllll}
\alpha^\star: {\P^3} & \to           & {\P^3}^\star    \\
X=(x_0,\ldots,x_3)\R & \mapsto & U            &= \R(u_0,u_1,u_2,u_3), \\
& &  &=\R(-(x_1^2+x_2^2+x_3^2),x_0x_1,x_0x_2,x_0x_3). \\
\end{array}
\ee
The base-points of $\alpha^\star$ are the origin $O=(1,0,0,0)\R$ and all points of the 
{\sf absolute conic} $j:x_0=0, x_1^2+x_2^2+x_3^2=0$, consisting of the circular points
at infinity. Considering the isotropic cone $\Delta$ as point set, we have 
$j=\Delta\cap\omega$.

For the practical examples in Section~\ref{addexam:sec} we use the fact that $\alpha$ 
and $\alpha^\star$ can be represented as composition of an inversion
and a polarity, see Figure~\ref{decompose_alpha_fig}. 
Let $\sigma:\P^3\to\P^3$ be the inversion at $S^2$, and let $\pi:{\P^3}^\star\to\P^3$ 
be the polarity with respect to $S^2$, and $\pi^\star: \P^3 \to {\P^3}^\star$ its dual
map. These maps satisfy $\sigma=\sigma^{-1}$ and $\pi\circ\pi^\star=\mbox{id}$,
and their coordinate representations are 
\be \label{polarity:eq}
\sigma: X=(x_0,\ldots,x_3)\R &\mapsto&
\sigma(X)=(x_1^2+x_2^2+x_3^2,x_0x_1,x_0x_2,x_0x_3)\R, \nonumber\\
\pi: U=\R(u_0,\ldots,u_3) &\mapsto& 
\pi(U) = X = (-u_0,u_1,u_2,u_3)\R,\\
\pi^\star: X=(x_0,\ldots,x_3)\R &\mapsto& 
\pi^\star(X) = U = \R(-x_0,x_1,x_2,x_3).\nonumber 
\ee
\begin{proposition}\label{footpoint_composition}
The foot-point map $\alpha$ and its inverse map $\alpha^\star$ from 
equations~\eqref{footpoint_map2:eq} and~\eqref{footpoint_map3:eq} are 
represented as compositions of the maps $\sigma, \pi$ and $\pi^\star$
from equation~\eqref{polarity:eq} by
\bes
\alpha = \sigma\circ\pi, \mbox{ and } \alpha^\star = \pi^\star\circ\sigma.
\ees
\end{proposition}

Given a rational unit vector field $\nv(u,v)\subset S^2$ and a rational function $r(u,v)$,
the map $\phi$ from~\eqref{offset:rep} creates a rational offset surface $F^\star$, 
and the map $\gamma$ from~\eqref{conch:rep} creates a rational conchoid surface $G$. 
Composing these maps with $\alpha$, we have $\alpha\circ\phi = \gamma$ 
and $\alpha^\star\circ\gamma = \phi$.
The relations between $\alpha$, $\phi$ and $\gamma$ are displayed in the
following diagrams.

\begin{equation}\label{footpoint_map}
\begin{diagram}
\node{} \node{S^2\times\R} \arrow{sw,l}{\phi} \arrow{se,l}{\gamma} \node{}\\
\node{\cal E} \arrow[2]{e,b}{\alpha} \node{} \node{\cal P} 
\end{diagram}
\quad\quad
\begin{diagram}
\node{} \node{S^2\times\R} \arrow{sw,l}{\phi} \arrow{se,l}{\gamma} \node{}\\
\node{\cal E}  \node{} \node{\cal P} \arrow[2]{w,b} {\alpha^\star} \node{} 
\end{diagram}
\end{equation}

\begin{theorem} \label{corr-theo:theo}
Let $F^\star\subset {\cal  E}$ be a dual surface, and let 
${o}^\star$ be the offset map~\eqref{offset_map:eq}, with $o^\star(F^\star)=F_d^\star$,
for arbitrary offset distance $d$. 
Let ${c}$ be the conchoidal map~\eqref{conchoid_map:eq}, with $c(G) = G_d$.
Then the family of surfaces $F_d^\star$ is mapped by $\alpha:{\cal  E} \to {\cal P}$ 
to a family of surfaces $G_d=\alpha(F_d^\star)$,
being conchoid surfaces of $G$ at distance $d$, with respect to the chosen reference
point $O$ of $\alpha$. Likewise, $\alpha^\star$ maps a family of conchoid surfaces $G_d$
to a family of offset surfaces $F_d^\star = \alpha^\star(G_d)$.
The following diagrams are commutative.
\begin{equation}\label{comm_diagram}
\begin{xy}
  \xymatrix{
   &  F^{\star}   \ar[r]^\alpha \ar[d]_{o^\star} &  G \ar[d]^c  \\
   &  F^{\star}_d \ar[r]_\alpha               &   G_d
  }
\end{xy}
\quad\quad
\begin{xy}
  \xymatrix{
   &  F^{\star} \ar[d]_{o^\star}  &  G \ar[d]^c \ar[l]^{\alpha^\star}  \\
   &  F^{\star}_d                 &   G_d       \ar[l]_{\alpha^\star} 
  }
\end{xy}
\end{equation}
\end{theorem}

\begin{proposition}\label{corr-theo2:theo} 
Let $F$ be an irreducible algebraic surface and $G=\alpha(F^\star)$. 
Then, the offset surfaces $F_d$ of $F$ are birationally equivalent to 
the conchoid surfaces $G_d$ of $G$.
\end{proposition}
This statement is implied by the birationality of the map $\alpha$. 
But other properties, as the degree, cannot be translated directly. 
There are relations between the degree of a dual surface $F^\star$ and the 
degree of the image surface $G=\alpha(F^\star)$, as outlined in the following.
Let $\ol{F}^\star$ be the zero-set of the irreducible polynomial~\eqref{fstar_implicit}.
Applying $\alpha$, which means inserting the representation~\eqref{footpoint_map3:eq}
into $\ol{F}^\star$ yields the polynomial 
\be\label{Gimplicit:eq}
\wt G(x_0,\xv) &= & (-\xv^2)^nf_0 + (-\xv^2)^{n-1}x_0f_1(\xv) + \ldots + 
(-\xv^2)^{k}x_0^{n-k}f_{n-k}(\xv)+ \ldots + \nonumber\\
& &(-\xv^2)x_0^{n-1}f_{n-1}(\xv) + x_0^nf_n(\xv).
\ee 
According to the fact that $\alpha$ possesses exceptional elements, the polynomial
$\wt G$ might have factors of the form $x_0^r$ or $(\xv^2)^k=(\xv\cdot\xv)^k$.
Since we are not interested in these factors, the pedal surface 
$\ol G = \alpha(\ol F^\star)$ is defined as the zero-set of the irreducible 
component of $\wt G$.
We analyze the degree of the image surface $\ol G$, depending on the position of
$\ol F^\star$ with respect to the exceptional planes of $\alpha$.
 
If $\omega$ is an $r$-fold plane of $\ol {F^\star}$, thus $f_0=\ldots f_{r-1}=0$,
but $f_r\not=0$, then the polynomial $\wt G$ contains the factor $x_0^r$.
If $f_n$ has the factor $\uv^2$, it also appears in $\wt G$. 
More generally, consider the case where the polynomials $f_n, \ldots, f_{n-k}$ have 
a common factor. In detail, let
\bes
f_{n-j}(\uv) = (\uv^2)^{k-j} h_{n-2k+j}(\uv), \mbox{ for } j=0, \ldots, k,
\ees   
with $\gcd(h_{n-2k+j}, \uv^2)=1$. Thus the polynomial~\eqref{Gimplicit:eq} 
has the factor $(\xv^2)^k$, and the irreducible component $G$ is of degree $2n-2k$.
Summarizing, the degree of $\ol G$ is $2n-r-2k$, where $r$ is the
multiplicity of $\omega$ and $k$ is the multiplicity of the cone
$u_0=0, \uv^2=0$ for $\overline{F}$.

In a similar way one can deduce the defining equation of a surface 
$\ol F^\star=\alpha^\star(\ol G)$, starting from a homogeneous polynomial
defining $\ol G$. Exchanging $F^\star$ by $G$ and dual coordinates $u_i$
by $x_i$, one obtains the degree of the dual object $\ol F^\star$. 
Hereby, one has to exchange also $\omega$ by $O$ and the isotropic cone $\Delta\ni O$
by the conic $j\subset\omega$.

\section{Elementary examples}\label{elemexam:sec}

To illustrate the relation between the offset and conchoid construction we
present two elementary examples. The first discusses conchoid surfaces $G_d$ 
of a plane $G$ with respect to a reference point $O\notin G$ and
the corresponding offset surfaces, which are offsets of a paraboloid
of revolution $F$ with $O$ as focal point.
The second example considers offsets of a sphere $F$ and
their corresponding conchoid surfaces,
which are typically rotational surfaces with a Pascal curve as
profile, and double point at $O$.

\subsection{Conchoid surfaces of a plane}

 Consider the plane $G: z=1$ and the reference point
 $O=(0,0,0)$. To compute the conchoid surfaces of $G$, one might use
 the trigonometric parameterization 
 \be\label{trigpara_s2:eq}
 \nv(u,v)=(\cos u \cos v, \cos v
 \sin u, \sin v )
 \ee  
 of $S^2$.
 Then $G$ and its conchoid surfaces $G_d$
 admit the trigonometric polar representations 
 \be\label{conch_plane:ex}
 \gv(u,v)&=&r(u,v)\nv(u,v), \mbox{ and }\\
 \gv_d(u,v)&=&(r(u,v)+d)\nv(u,v), \mbox{ with }  r(u,v) =\frac{1}{\sin v}. 
 \nonumber
 \ee
 Applying the map $\alpha^\star$ from~\eqref{footpoint_map_inv:eq} 
 to these parameterizations gives dual parametric representations of 
 $F$ and $F_d$. The tangent planes of these surfaces read
 \bes
 E(u,v): \xv\cdot\nv(u,v) = r(u,v), \mbox{ and } E_d(u,v):
 \xv\cdot\nv(u,v) = r(u,v)+d.
 \ees 
 The surface $F: x^2+y^2+4z=4$ is a paraboloid of revolution with 
 focal point $O$ and axis $z$, 
 and $F_d$ are its offset surfaces. To illustrate the projective
 version of $\alpha$, we turn to the projective setting.
 Let $X=(x_0,\ldots,x_3)$ and $U=(u_0,\ldots, u_3)$, 
 the surfaces $\ol G\subset \P^3$ and $\ol F^\star\subset{\P^3}^\star$ 
 are the zero-sets of the homogeneous polynomials 
 \bes
 \ol{G}(X)=(x_3-x_0) 
 \mapsto 
 \alpha^\star(\ol G) = F^\star(U)= u_1^2+u_2^2+u_3^2 + u_0u_3,
 \ees
 To obtain the implicit equation of
 the conchoid surface $\ol G_d$, one has to eliminate the parameters
 $u$ and $v$ from~\eqref{conch_plane:ex}. One obtains the polynomial 
 of degree four, 
 \bes
 \ol G_d(X) = d^2x_0^2x_3^2 - (x_1^2+x_2^2+x_3^2)(x_0-x_3)^2.
 \ees
 Since the highest power of $x_0$ is two, the origin is a double 
 point of $\ol G_d$. This tells us that $\alpha^\star(\ol G_d)$
 contains the factor $u_0^2$. Additionally we may divide by $(u_1^2+u_2^2+u_3^2)$, 
 such that finally the offset $\ol{F}_d^\star$ of the paraboloid
 $\ol F^\star$ is the zero-set of the polynomial
 \bes
 \ol{F}_d^\star(U) = d^2u_3^2(u_1^2+u_2^2+u_3^2) - (u_1^2+u_2^2+u_3^2 + u_0u_3)^2.
 \ees
 Since the highest power of $u_0$ is two, the plane at infinity $\omega=\R(1,0,0,0)$
 is a two-fold tangent plane of $\ol{F}_d^\star$.
 Figure~\ref{parabola_offset_fig} provides an illustration of that
 example in the 2d-case.

\subsection{Offsets of a sphere and corresponding conchoid surfaces}
\label{sphere_offsets}

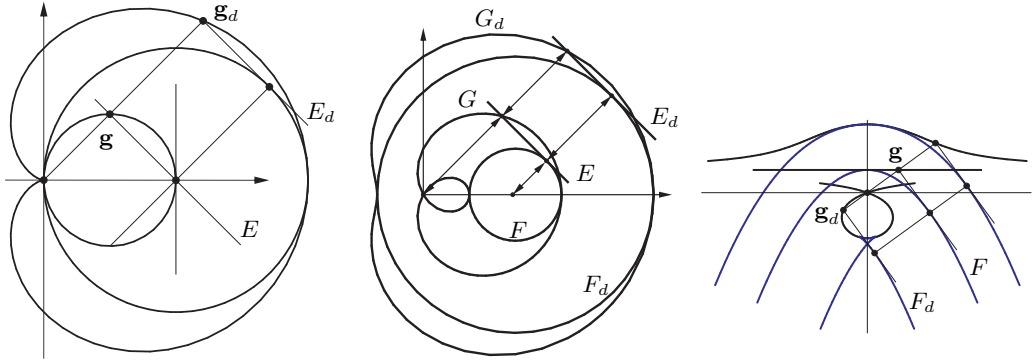
\begin{figure}[ht]
 \subfigure[Offsets of a pencil and corresponding conchoid
 curves]{\label{linepencil_fig}
 \begin{overpic}[width=.3\textwidth]{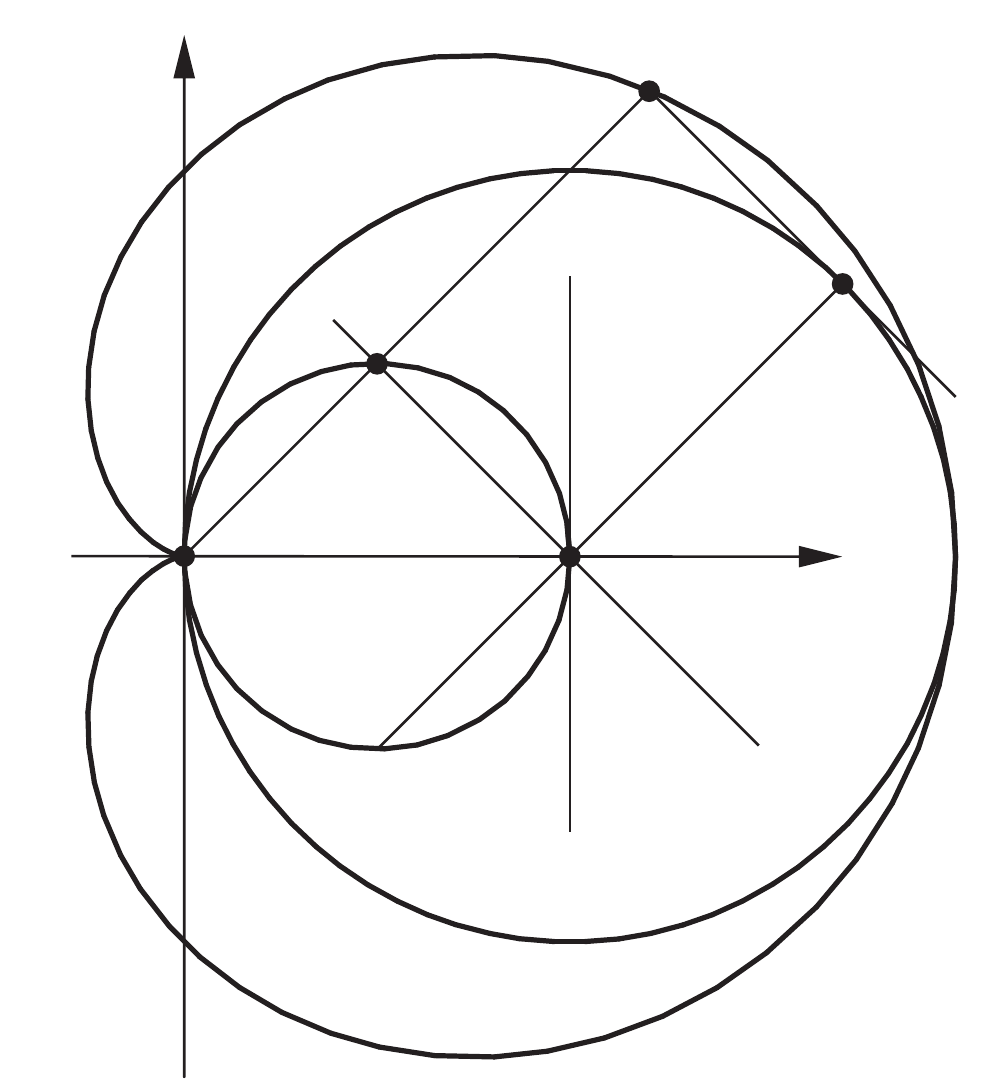}%
 {\small
 \put(32,58){$\gv$} \put(70,33){$E$}  \put(88,65){$E_d$}
 \put(63,93){$\gv_d$}
 }
 \end{overpic}
 }
 \subfigure[Offsets of a circle and corresponding conchoid
 curves]{\label{circleoffset_fig}
 \begin{overpic}[width=.3\textwidth]{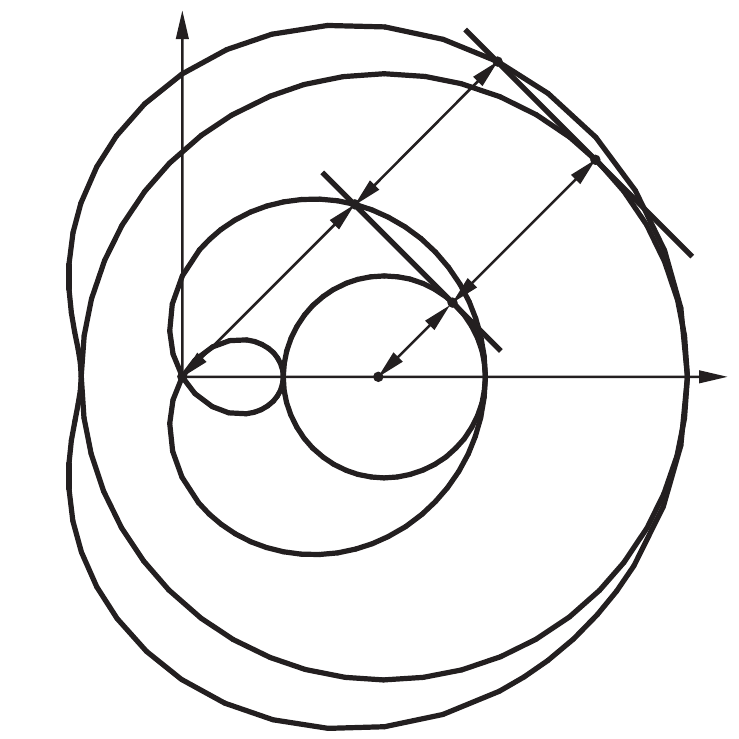}%
 {\small
 \put(50,36){$F$} \put(72,20){$F_d$}
 \put(35,75){$G$} \put(70,53){$E$}  \put(92,70){$E_d$} \put(40,99){$G_d$}
 }
 \end{overpic}
 }
 \subfigure[Offsets of a parabola and conchoids of a
 line]{\label{parabola_offset_fig}
 \begin{overpic}[width=0.3\textwidth]{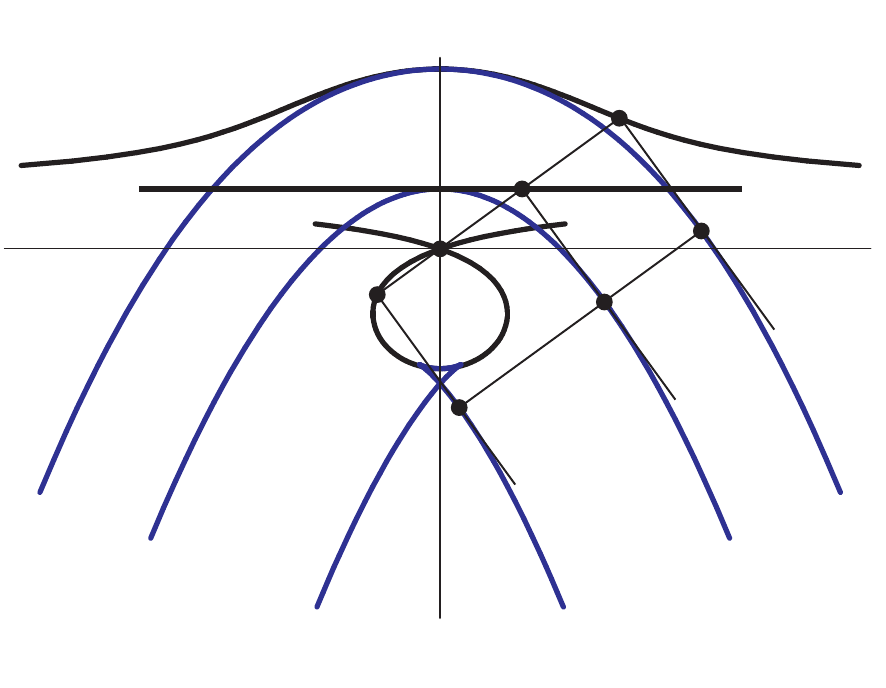}%
 {\small
 \put(62,15){$F_d$} \put(80,25){$F$}
 \put(34,42){$\gv_d$} \put(56,60){$\gv$}  }
 \end{overpic}
 }
 \caption{Relation between offsets and conchoids}
 \label{geom_prop:fig}
 \end{figure}

 Consider the sphere $F:(x-m)^2+y^{2}+z^2-R^2=0$. Its offset surfaces
 $F_d$ are spheres
 with same center and Radius $R+d$. For a 2d-illustration see
 Figure~\ref{circleoffset_fig}.
 To establish the correspondence with a family
 of conchoid surfaces $G_d$, we have to consider the dual surfaces
 $F_d^\star$, determined by the polynomial
 \bes
 \ol F_d^\star(u_0,\uv) &=& ((R+d)^2-m^2)u_1^2 + (R+d)^2(u_2^2 +
 u_3^2) - 2mu_0u_1 - u_0^2.
 \ees
 Substituting $d=0$ gives $\ol F^\star$.
 The bi-rational map $\alpha$ maps surfaces $\ol F_d^\star$
 to the family of conchoid surfaces $\ol G_d$, which are the zero-set
 of the polynomial of degree four,
 \bes
 \ol G_d(x_0, \xv) &=& x_0^2(x_1^2((R+d)^2-m^2) + (R+d)^2(x_2^2+x_3^2))
  + 2mx_0x_1(\xv^2) - (\xv^2)^2.
 \ees
 Since the highest power of $x_0$ is two, the origin $O=(1,0,0,0)\R$
 is a double point of $\ol G_d$. Substituting $d=0$ gives $\ol G$. 
 We notice that when letting $R=0$ and $d=0$, the sphere $F$
 degenerates to a (double traced) bundle of planes passing through the point $M=(m,0,0)$. 
 The corresponding 2d-example is illustrated in Figure~\ref{linepencil_fig}.
 Its defining polynomial is $\ol F^\star(u_0,\uv)=(u_0+mu_1)$.
 The corresponding surface $\ol G=\alpha(\ol F^\star)$ is the zero-set
 of the polynomial $G(x_0,\xv)=\xv\cdot\xv-mx_0x_1$. 
 Thus, $\ol G$ is a sphere with diameter $OM$.

Dual rational parameterizations of $F_d$ and rational polar representations 
of $G_d$ may be derived as follows. Consider again the trigonometric
parameterization~\eqref{trigpara_s2:eq} of $S^2$.
Tangent planes $E_d$ of $F_d$ and a parameterization $\gv_d$ are consequently 
given by
\bes
E_d(u,v):& & \xv\cdot\nv(u,v) = r(u,v)+d,\\
\gv_d(u,v)&=&(r(u,v)+d)\nv(u,v), \mbox{ with } 
r(u,v)=m\cos u \cos v + R.
\ees

\section{Pedal surfaces and inverse pedal surfaces}\label{addexam:sec}

Consider a surface $F$, its offsets $F_d$ and a fixed reference point $O=(0,0,0)$.
The surface $G_d=\alpha(F_d^\star)$ is called {\sf pedal surface} of $F_d$,
and consists of the foot-points at the tangent planes of $F_d$ with respect to $O$.
According to Theorem~\ref{corr-theo:theo}, the surfaces $G_d$ are the conchoid surfaces 
of $G$ with respect to $O$.
Conversely, the surface $F_d^\star=\alpha^\star(G_d)$ is called {\sf negative pedal surface}
of $G_d$ with respect to $O$. Since the operation applied is rather the inverse than 
the negative, we use the notion {\sf inverse pedal surface}. 
Consider points $X\in G_d$, then the tangent planes
of $F_d$ are perpendicular to $OX$ and contain $X$.

There exist several interesting examples, both for rational offset
surfaces and for rational conchoid surfaces and their counterparts with respect to
$\alpha$ and $\alpha^\star$. These maps apply to translate geometric properties 
between offset surfaces and conchoid surfaces and vice versa, as stated 
in Proposition~\ref{corr-theo2:theo}. The following subsections discuss 
pedal surfaces and inverse pedal surfaces of ruled surfaces and quadrics.

\subsection{Pedal surfaces of rational ruled surfaces}
\label{pedal_rs:sec}

\begin{figure}[ht]
\centering
\subfigure[Circle as pedal curve of a pencil of planes]{\label{pedal_planepencil}
\begin{overpic}[width=.2\textwidth]{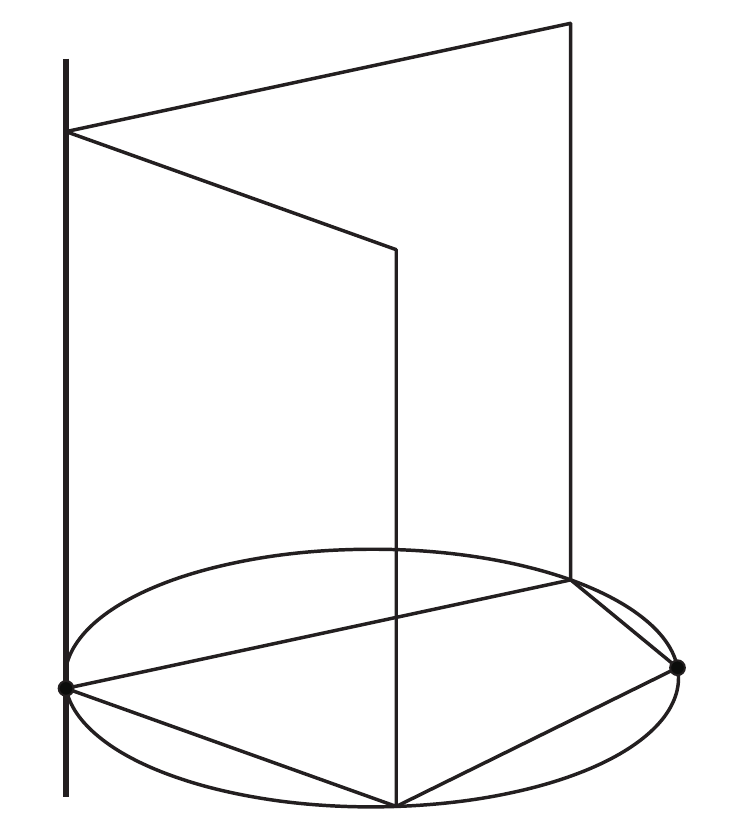}%
{\small  \put(-2,15){$D$}  \put(83,15){$O$} \put(72,4){$a$}}
\end{overpic}
}
\hfil
\subfigure[Parabola as inverse pedal curve of a line]{\label{invpedal_line}
\begin{overpic}[width=.3\textwidth]{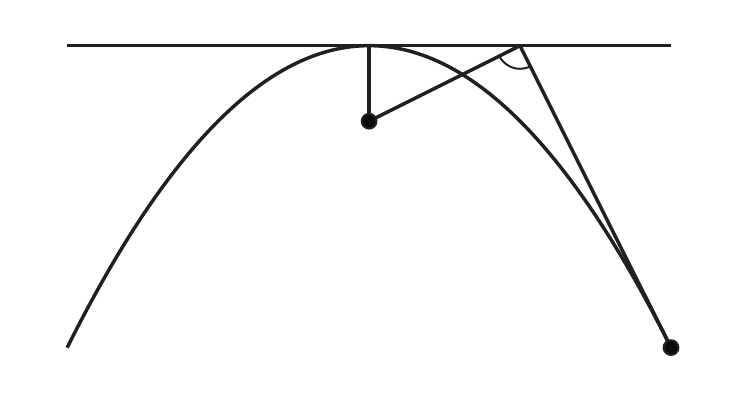}%
{\small    \put(42,30){$O$} \put(21,49){$l$} \put(71,49){$X$}}
\end{overpic}
}
\caption{Pedal and inverse pedal curve of a line}
\label{pedal_invpedal_line:fig}
\end{figure}
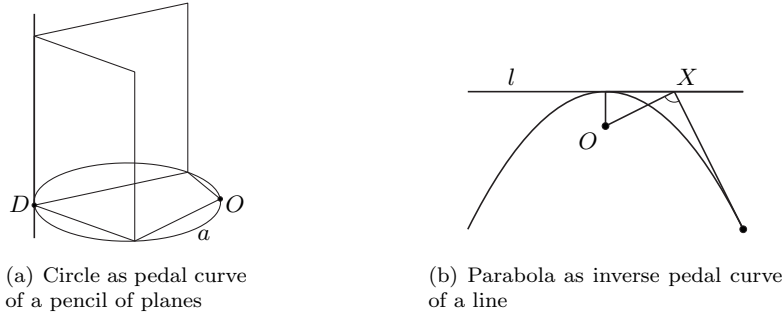

Consider a skew ruled surface $\ol F\subset\P^3$. Its tangent planes form 
the dual surface $\ol F^\star$ which is itself a skew ruled surface in ${\P^3}^\star$. 
Consider a fixed line $l\in F$, and let $D$ be the foot-point of $O$ on $l$.
Actually we interpret $l$ as pencil of planes, and not as point set.
According to Thales' theorem, the pedal curve $\alpha(l)$ is a circle with 
diameter $OD$, in a plane perpendicular to $l$, see Figure~\ref{pedal_planepencil}.
Consequently, 
the pedal surface $\ol G=\alpha(\ol F^\star)$ contains a one-parameter
family of circles in planes perpendicular to the lines $l$ of $F$.

An implicit representation of a surface $\ol F$ does not tell us if the surface
is ruled or not, except for some special cases. 
Thus we derive a construction for $\ol G=\alpha(\ol F^\star)$ in terms 
of a parameterization of $\ol F^\star$. 
It is known, see for instance~\cite{pet1}, that a rational ruled surface is 
a rational offset surface in the sense of Definition~\ref{ratoffset:def}.
But typically it makes some effort to construct a rational offset parameterization 
for $\ol F^\star$. Thus we start to demonstrate the construction in terms of 
a general parameterization of a ruled surface.

Let $\fv(u,v)=\cv(u) + v\ev(u)$ be an affine parameterization of $F$. The foot-points $\dv$
on the generating lines of $F$ with respect to $O$ have to satisfy $\dv\cdot\ev=0$.
This implies that $v = -(\cv\cdot\ev)/(\ev\cdot\ev)$ is the parameter value of $\dv$ 
and we have 
\be\label{rsfootpointcrv:eq}
\dv(u) = \cv(u) - \frac{\cv(u)\cdot\ev(u)}{\ev(u)\cdot\ev(u)} \ev(u).
\ee
We exchange the directrix curve $\cv$ by $\dv$. The circles $\alpha(l(u))$ have centers 
at $\mv(u)=\frac{1}{2} \dv(u)$, and their carrier planes $\eps(u)$ are orthogonal to $l$, 
thus $\eps(u):\xv\cdot\ev(u) = 0$.

To find a parameterization of the pedal surface $G$, we need to parameterize the
family of circles $\alpha(l(u))$. For that we need an orthonormal frame in $\eps(u)$.
Consider two orthonormal vectors $\av(u), \bv(u)$, such that 
$\{ \ev / \| \ev \|, \av,\bv \}$ is an orthonormal basis in $\R^3$. 
Obviously these vectors are solutions of $\ev(u)\cdot\xv=0$.
Since the radius of $\alpha(l(u))$ is given by $r(u) = \|\dv \|/2$, 
an affine parameterization of $G$ reads
\bes
\gv(u,v) = \frac{1}{2}\left( \dv(u) + \|\dv(u) \| \av(u) \cos v + 
                                      \|\dv(u) \| \bv(u) \sin v \right).
\ees
Typically this is a non-rational parameterization even for 
rational input surfaces, since the radius $r(u)$ and also $\av(u)$ and $\bv(u)$ 
are typically non-rational.
But on the other hand, $G$ is a rational surface, since it carries a 
rational family of circles given by the equations
\bes 
\|\xv-\mv(u)\|^2 - r(u)^2=0, \mbox{ and } \xv\cdot\ev(u) = 0.
\ees
Rational parameterizations of a rational one-parameter family of conics
can be computed explicitly, see for instance~\cite{schicho1}. Nevertheless, 
to represent the whole family of conchoid surfaces $G_d$ by a rational 
parameterization, we have to start with a representation of 
the tangent planes of $F$, involving a rational unit normal vector. 
Computing the partial derivatives of $\fv(u,v)=\cv(u)+v\ev(u)$
gives $\fv_u = \dot\cv(u) + v\dot\ev(u)$ and $\fv_v(u,v) = \ev(u)$. 
Thus the normal $\nv(u,v)$ reads
\bes
\nv(u,v) = \nv_1(u) + v\nv_2(u), \mbox{ with }
\nv_1(u) = \dot\cv(u)\times\ev(u), \mbox{ and }
\nv_2(u) = \dot\ev(u)\times\ev(u).
\ees
In case that $\nv_1(u)$ and $\nv_2(u)$ are linearly dependent for all $u$,
$F$ is a developable ruled surface, which can be considered as the envelope of its 
one-parameter family of tangent planes. In this case, the pedal 'surface' $\alpha(F^\star)$
degenerates to a curve, and therefore we exclude this case. 

Let $\nv_1$ and $\nv_2$ be linearly independent, which is equivalent to  
$\det(\dot\cv, \ev, \dot\ev)\not=0$, since $\ev(u)$ is assumed to be a non-zero direction
vector field. In order to represent the family of offset surfaces $F_d$ by a rational 
parameterization, we have to construct a rational unit normal vector field of $F$. 
To obtain this, we require $\|\nv(u,v)\|^2 = w^2$, with some function $w(u,v)$ 
to be determined.
To simplify the practical computation, we exchange the directrix curve $\cv(u)$ by the
striction curve $\sv(u)=\cv(u)+v_s(u)\ev(u)$, having the property that the normal vector 
$\nv_s$ at $\sv$ is orthogonal to $\nv_2=\dot\ev\times\ev$. 
This condition determines the parameter 
$$
v_s(u)=-\frac{(\dot\cv\times\ev)\cdot(\dot\ev\times\ev)}{(\dot\ev\times\ev)^2}.
$$
Representing $F$ with $\sv$ as directrix curve by $\fv(u,v)=\sv(u)+v\ev(u)$ 
yields that $\nv=\nv_s + v\nv_2$, with $\nv_s\perp\nv_2$. 
The condition $\| \nv \|^2 = w^2$ turns into $\nv_s^2 + v^2\nv_2^2 = w^2$, where
we consider $v$ and $w$ as affine coordinates in $\R^2$. By the substitution $v=y_2/y_1$ and 
$w =y_0/y_1$, this quadratic equation becomes
\be\label{conicfam:eq}
a(u): a_1(u)y_1^2 + a_2(u)y_2^2 - y_0^2=0, \mbox{ with } a_1(u)=\nv_s(u)^2, \mbox{ and }
a_2(u) = \nv_2(u)^2.
\ee
Equation~\eqref{conicfam:eq} defines a {\sf real rational} family of conics in the projective
plane $\P^2$, with coordinates $y_0,y_1,y_2$. There exists a parameterization 
$\yv(u,t)=(y_0,y_1,y_2)(u,t)$, satisfying equation~\eqref{conicfam:eq} identically, 
in a way that $\yv(u_0,t)$ is a parameterization
of the conic $a(u_0)$, for any fixed $u_0\in\R$, see for instance~\cite{schicho1}.   
Consequently one has constructed a parameterization 
\be\label{ratoffset_rs:eq}
\fv(u,t) = \sv(u) + \frac{y_2(u,t)}{y_1(u,t)}\ev(u), \mbox{ with } 
\|\nv(u,t) \|=\frac{y_0(u,t)}{y_1(u,t)}.
\ee 
Now, the norm of the normal vector $\nv(u,t)$ is a rational function in the surface parameters
$u$ and $t$. We represent the ruled surface $F$ and its offset surfaces $F_d$ as envelopes of 
their tangent planes 
\bes
E(u,t)&:& \xv\cdot\nv(u,t) = \fv(u,t)\cdot\nv(u,t), \nonumber\\
E_d(u,t)&:&  \xv\cdot\nv(u,t) = \fv(u,t)\cdot\nv(u,t) + d\|\nv(u,t) \|.
\ees 
In other words, $E_d = \R(-\fv\cdot\nv-d\|\nv\|, \nv)$ is a dual rational parameterization
of the family of offset surfaces $F_d^\star$ of the given ruled surface $F$, and analogous
to equation~\eqref{offset_surf_dual2:eq}, up to the normalization of $\nv$.
Applying the foot-point map $\alpha$ 
yields a rational parameterization of the corresponding family of rational conchoid 
surfaces $G_d=\alpha(F_d^\star)$. In terms of homogeneous coordinates we have 
\bes
E_d = \R(-\fv\cdot\nv-d\|\nv\|, \nv) \mapsto
\alpha(E_d) &=& (\nv^2, (\fv\cdot\nv+d\|\nv \|)\nv)\R.
\ees 
The according affine rational polar representation of $G_d$, with $\fv(u,t)$ 
from equation~\eqref{ratoffset_rs:eq} and $\nv(u,t)$ its corresponding normal
vector field, reads
\bes
\gv_d(u,t) = \left(\frac{\fv\cdot\nv}{\| \nv \|} + d\right)\frac{\nv}{\|\nv \|}.
\ees 

\begin{proposition}
The pedal surfaces $G_d$ of the family of offset surfaces $F_d$ of 
non-developable rational ruled surfaces $F$ admit rational polar representations 
with respect to 
any chosen reference point $O$. The pedal surface $G$ of $F$ is generated by
a rational family of circles in planes through $O$ and perpendicular to 
$F$'s generating lines $l$. The conchoid surface $G_d$ of $G$ is generated by
the planar conchoid curves of that circles lying on $G$. 
\end{proposition}

\subsection{Inverse pedal surfaces of rational ruled surfaces}
\label{invpedal_rs:sec}

Consider a rational ruled surface $G$ as set of points, and a fixed reference point 
$O=(0,0,0)$. Applying $\alpha^\star$ yields
a two parameter family of tangent planes $F^\star=\alpha^\star(G)$,
whose envelope $F$ is the inverse pedal surface of $G$. While the ruled
surface $F$ in Section~\ref{pedal_rs:sec} has been assumed to be non-developable,
since it has been considered as two-parameter family of tangent planes,
the ruled surface $G$ in this current section may be a developable surface as well,
since it is considered as two-parameter family of points. Thus, $G$ might be a 
tangent surface of a curve, a cylinder, a cone or even a plane. The latter
case makes sense only if $O\notin G$.

Let $l\subset G$ be a generating line, and consider points $X\in l$.
The one-parameter family of planes $\alpha^\star(X)$ through $X$, whose normal vector
is $OX$, envelope a parabolic cylinder $P$, with $O$ as focal point of the cross section
parabola $p$ in the plane connecting $O$ and $l$, see Figure~\ref{invpedal_line}.  
The vertex of the parabola $p$ is obviously the foot-point of $O$ on $l$. 
Consequently, the inverse pedal surface $F$
of a ruled surface $G$ is the envelope of a one parameter family of parabolic cylinders.
The orthogonal cross sections of these parabolic cylinders with planes through $O$ have 
$O$ as common focal point.

Consider a general rational ruled surface parameterization $\gv(u,v)=\dv(u)+v\ev(u)$ 
of $G$, where $\dv(u)$ is the foot-point curve with respect to $O$, compare 
equation~\eqref{rsfootpointcrv:eq}. The vertices of above mentioned parabolas are $\dv(u)$.
The tangent planes of $F$ are $E(u,v): (\xv-\gv)\cdot\gv=0$, but since $\gv$ is typically 
not a rational polar representation, the normal vector of $E(u,v)$ does not have rational norm.
In detail we have
\be\label{invpedalrs:eq}
E(u,v)&:& \xv\cdot(\dv(u) + v \ev(u)) = \dv(u)^2 + v^2\ev(u)^2, \nonumber\\
E_u(u,v)&:& \xv\cdot(\dot\dv(u) + v\dot\ev(u)) = 
2(\dv(u)\cdot\dot\dv(u) + v^2\ev(u)\cdot\dot\ev(u)), \nonumber\\
E_v(u,v)&:& \xv\cdot\ev(u) = 2v\ev(u)^2.
\ee  
The solution $\fv(u,v)$ of~\eqref{invpedalrs:eq} is a rational representation 
of the inverse pedal surface $F$ of the ruled surface $G$. Although $F$ is a rational
offset surface in the sense of Definition~\ref{ratoffset:def}, the unit normal 
vector field $\nv(u,v)$, corresponding to $\fv(u,v)$ is typically non-rational. 
Nevertheless, the intersection $E\cap E_v$ gives the generating lines of the 
one-parameter family of parabolic cylinders $P(u)$.
These cylinders $P(u)$ admit the parameterization 
\bes
\qv(u, v, \lambda) = \pv(u,v) + \lambda \av(u), 
\ees 
with cross-section parabolas $\pv$ and $\av$ as direction vectors of its generating lines, 
given by
\bes
\pv(u,v) = \left( 1-\frac{v^2\ev^2}{\dv^2} \right) \dv + 2v\ev, 
\mbox{ and } \av(u) = \dv(u)\times\ev(u).
\ees
The intersection of $\qv(u,v,\lambda)$ with planes $E_u(u,v)$ from~\eqref{invpedalrs:eq}
determines $\lambda$ as a rational function in $u$ and $v$, and finally gives the 
parameterization $\fv(u,v)$ of $F$. The $v$-lines of $\fv(u,v)$, that are 
the characteristic curves of the family of cylinders $P(u)$,
are rational curves of degree at most three. 

The development of a parameterization $\fv(u,v)$ of $F$ with the property that 
its unit normal vectors are rational, is more involved. To construct such a 
parameterization, one needs to start with a rational parameterization 
$\gv(u,v)$ of $G$, whose norm $\|\gv(u,v)\|$ is a rational function in $u$ and $v$.
We do not give the detailed construction here but refer to~\cite{pet10}. To outline
the method, one studies the squared norm 
\be\label{normrs:eq}
\|\gv(u,v) \|^2 = \dv(u)^2 + v^2\ev(u)^2 = w^2,
\ee
which shall be the square of a rational function $w(u,v)$. Substituting $v=y_1/y_2$
and $w=y_0/y_2$, equation~\eqref{normrs:eq} defines a rational one-parameter family of conics
\be\label{conicfam_invpedal:eq}
a(u): y_2^2\dv(u)^2 + y_1^2\ev(u)^2 - y_0^2=0.
\ee
As mentioned in Section~\ref{pedal_rs:sec} and proved in~\cite{schicho1}, these objects 
have rational parameterizations. More precisely, one can construct rational functions
$y_i(u,t)$, satisfying~\eqref{conicfam_invpedal:eq} identically, so that the $t$-lines
of $(y_0,y_1,y_2)(u,t)$ are the conics $a(u)$. Considering the reparameterization
$v=y_1/y_2$, one obtains the parameterization
\bes 
\gv(u,t) = \dv(u) + \frac{y_1(u,t)}{y_2(u,t)}\ev(u), \mbox{ with }
\|\gv(u,t)\| = \frac{y_0(u,t)}{y_2(u,t)}.
\ees 
Using this specific parameterization of $G$ instead of $\gv=\dv+v\ev$, the family of 
tangent planes $E(u,t): (\xv-\gv(u,t))\cdot\gv(u,t)=0$ possesses a normal
vector $\gv(u,t)$ with rational length $\|\gv(u,t)\|$. Consequently the intersection point 
$\fv(u,t)$ of the planes $E(u,t)\cap E_u(u,t)\cap E_t(u,t)$ is a rational parameterization
of $F$, whose corresponding unit normal vector $\nv(u,t)=\gv(u,t)/\|\gv(u,t)\|$ 
is rational.  

\begin{proposition}
The inverse pedal surfaces $F = \alpha^\star(G)$ of rational ruled surfaces $G$ 
admit rational parameterizations $\fv(u,t)$, whose corresponding unit normal vector 
$\nv(t,u)$ is rational, too.
The family of rational conchoid surfaces $G_d$ of the rational ruled surface $G$ 
is mapped by $\alpha^\star$ to a family of offset surfaces $F_d$ of $F$. While $F$ 
is the envelope of a rational one-parameter family of parabolic cylinders $P(u)$ 
with cross section parabolas $p(u)$, the offset surfaces $F_d$ are the envelopes 
of the offset cylinders $P_d(u)$, whose cross sections are offsets of the parabolas $p(u)$. 
\end{proposition}

\begin{figure}[ht]
\subfigure[Pl\"ucker's conoid $F$]{\label{pluecker_conoid}
\begin{overpic}[width=.3\textwidth]{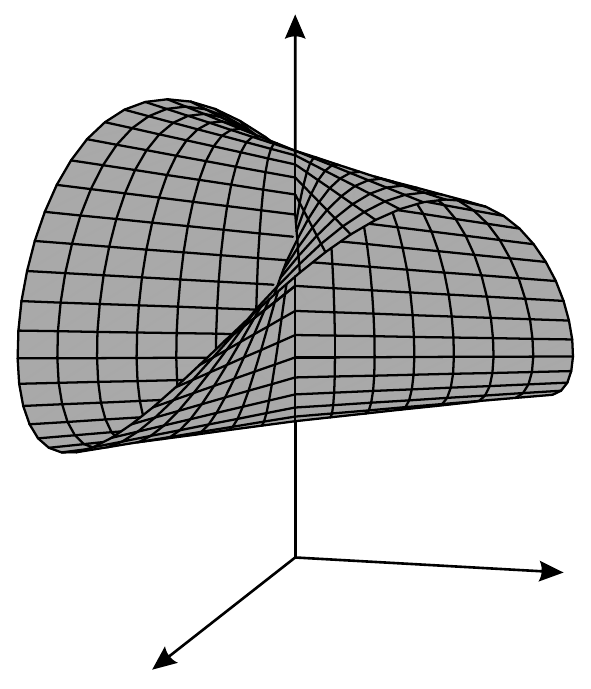}%
{\small  \put(65,75){$F$}  \put(42,12){$O$}}
\end{overpic}
}
\subfigure[Pedal surface of $F$]{\label{pedal_rs}
\begin{overpic}[width=.3\textwidth]{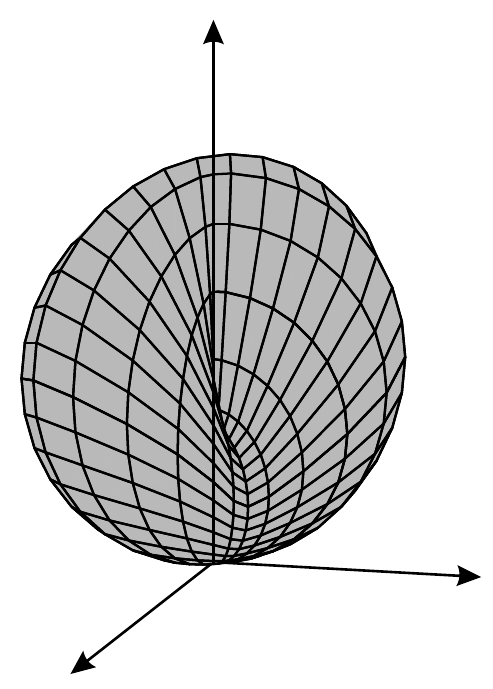}%
{\small    }
\end{overpic}
}
\subfigure[Inverse pedal surface of $F$]{\label{invpedal_rs}
\begin{overpic}[width=.3\textwidth]{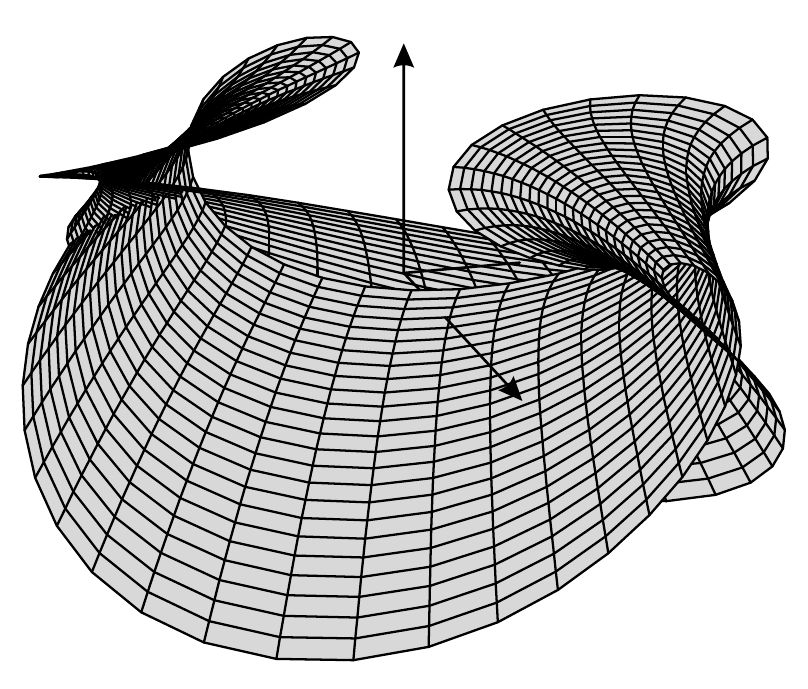}%
{\small  }
\end{overpic}
}
\caption{Pedal surface and inverse pedal surface of Pl\"ucker's conoid $F$ with 
respect to the origin $O$}
\label{ruled_surf:fig}
\end{figure}

\paragraph{Example}
We consider Pl\"ucker's conoid $F$, which is a ruled surface of
degree three, projectively equivalent to Whitney's umbrella. 
A possible parametrization of $F$ reads
\be\label{pluecker_para:eq}
\fv(r,u) = (r\cos u, r\sin u, \sin 2u) = \cv(u)+r\ev(u), 
\ee  
with directrix curve $\cv(u)=(0,0,\sin 2u)$ and direction
vectors $\ev(u)=(\cos u, \sin u,0)$. 
Considering $X=(x_0,\ldots,x_3)$ as homogeneous coordinates in $\P^3$,
the surface $F$ is the zero-set of the polynomial
\be\label{pluecker_eq:eq}
\ol F(X) = x_3(x_1^2+x_2^2)-2x_0x_1x_2.
\ee
We observe that the line $x_1=x_2=0$ is a double line of $\ol F$. Note that 
for the computation of the equations and parameterizations we choose $O$ as
symmetry point on the $z$-axis, whereas the illustrations in Figure~\ref{ruled_surf:fig}
have been generated for a diffreent reference point on the $z$-axis. The symmetric
position implying~\eqref{pluecker_eq:eq} gives more compact equations for offsets,
pedal surfaces and inverse pedal surfaces of $F$. 
 
To compute the pedal surface $G$ of $F$, we have to compute the tangent planes 
and the dual representation of $F$. From~\eqref{pluecker_para:eq} it follows that 
the tangent planes $E$ and $E_d$ of $F$ and its offset surface $F_d$ are represented by the homogeneous plane coordinates
\be\label{pluecker_para_dual:eq}
E(r,u) = \R(r\sin 2u, -2\sin u\cos 2u, 2\cos u\cos 2u, -r),\nonumber\\
E_d(r,u) = \R(r\sin 2u - d\sqrt{\nv\cdot\nv}, 
-2\sin u\cos 2u, 2\cos u\cos 2u, -r), 
\ee
where $\nv=(-2\sin u\cos 2u, 2\cos u\cos 2u, -r)$. Note that $E(r,u)$ is not a rational
offset parameterization of $F$. 
Either by eliminating $\lambda, r$ and $u$ from the equations $\lambda u_i=E_i$, 
and $\lambda u_i={E_d}_i$, $i=0,\ldots,3$, or by applying the method outlined 
at the end of Section~\ref{dualrep_offset:sec}, one obtains implicit dual 
representations of $F$ and its offset surface $F_d$ by
\bes
\ol F^\star(U) &=& u_0(u_1^2 + u_2^2) - 2u_1u_2u_3, 
\mbox{ with } U=(u_0,\ldots, u_3), \nonumber\\
\ol F_d^\star(U) &=& d^2(u_1^2 + u_2^2)^2(u_1^2+u_2^2+u_3^2) 
- (u_0(u_1^2+u_2^2)-2u_1u_2u_3)^2 \nonumber \\
&=& d^2(u_1^2 + u_2^2)^2(u_1^2+u_2^2+u_3^2) - \ol F^\star(U)^2.
\ees
This tells that Pl\"ucker's conoid is not only of degree three but also of class three,
and its offsets are of class six.
To represent $F$ by a parameterization so that its normal vector has rational norm,
one considers the equation $\|\nv(r,u) \|^2=4\cos^2 2u + r^2 = w^2$. Its right hand side
is a family of conics in the plane $\R^2$ with coordinates $r$ and $w$. These conics are 
parameterized by  
\bes 
r(u,t) = \frac{2\cos 2u \cos t}{\sin t}, \mbox{ and }
w(u,t) = \|\nv \|=\frac{2\cos 2u}{\sin t}.
\ees 
Thus, the surface $F$ admits a parameterization $\fv(r(u,t), u)$, whose unit normal
vector is $\nv(u,t) = (-\sin u\sin t, \cos u\sin t, \cos t)$.
That reparameterization implies that the tangent planes~\eqref{pluecker_para_dual:eq}
of $F$ and its offsets $F_d$ are represented by
\bes
E_d(u,t) = \R(-\cos t\sin 2u - d, -\sin u\sin t, \cos u\sin t, \cos t).
\ees
With help of the transformation $\alpha$ and its inverse $\alpha^\star$, the defining 
polynomials of the pedal surfaces $\ol G$ and $\ol G_d$ follow by
\bes
\ol G(X) = 2 x_0 x_1 x_2 x_3 + (x_1^2+x_2^2)(x_1^2+x_2^2+x_3^2),\nonumber\\
\ol G_d(X) = d^2x_0^2(x_1^2+x_2^2)^2(x_1^2+x_2^2+x_3^2) - \ol G(X)^2.
\ees 
Additionally we want to give the implicit representations of the conchoid
surfaces of Pl\"ucker's conoid and their inverse pedal surfaces. Since Pl\"ucker's
conoid and its conchoid surfaces are now considered as point sets, we use a different 
notation as above. Let Pl\"ucker's conoid $A$ be the zero-set of the polynomial
\bes 
\ol A(X) = x_3(x_1^2+x_2^2)-2x_0x_1x_2,
\ees
which is just the same as equation~\eqref{pluecker_eq:eq}. Since $A$ is a rational
conchoid surface, there exists an affine parameterization $\av(u,v)$ with rational norm.
This and the parameterization of its conchoid surfaces read
\be\label{pluecker_conchoid_para:eq}
\av(u,v) &=& \left( \frac{2 \sin u \cos v^2 \sin v}{\cos(u)}, 
\frac{2\sin u \cos v \sin v^2}{\cos u}, 
2\cos v \sin v \right), \nonumber \\
\av_d(u,v) &=& \av(u,v) + d\frac{\av(u,v)}{\|\av(u,v) \|},  
\mbox{ with } \|\av \| = \frac{2\cos v\sin v}{\cos u}
\ee
The conchoid surfaces $A_d$ are of degree eight, and are the zero set 
of the polynomial
\bes
\ol A_d(X) &=& d^2 (x_1^2+x_2^2)^2 x_0^2x_3^2-
(x_1^2+x_2^2+x_3^2)(x_3(x_1^2+x_2^2)-2 x_0 x_1 x_2)^2, \nonumber\\
&=& d^2 (x_1^2+x_2^2)^2 x_0^2x_3^2-(x_1^2+x_2^2+x_3^2) \ol A(X)^2.
\ees
The inverse foot-point map $\alpha^\star$ transforms Pl\"ucker's conoid $\ol A$ 
and its conchoid surfaces $A_d$ to $\ol B^\star$ and a family of rational offset 
surfaces  $\ol B_d^\star$.
Rational offset parameterizations of $B^\star$ and $B_d^\star$ can be 
derived from the parameterization~\eqref{pluecker_conchoid_para:eq}.
The defining polynomials of these surfaces in dual coordinates read
\bes
\ol B^\star(U) = u_0u_3(u_1^2 + u_2^2) + 2 u_1u_2(u_1^2 + u_2^2 + u_3^2), \nonumber\\
\ol B_d^\star(U) = d^2(u_3^2)(u_1^2+u_2^2)^2(u_1^2+u_2^2+u_3^2) - {\ol B^\star}^2. 
\ees

\subsection{Pedal surfaces of quadrics}
   
Consider a quadric $\ol F\subset \P^3$ as family of tangent planes.
Quadratic cylinders and cones are excluded since their pedal 'surfaces'
are curves. 
The dual surface $\ol F^\star\subset {\P^3}^\star$ is also a quadric. Using
homogeneous plane coordinates $U=(u_0,u_1,u_2,u_3)$, a dual quadric 
is the zero-set of a quadratic polynomial 
\be\label{dualquadric:eq}
\ol F^\star(U) &=& U^T \cdot A \cdot U, \mbox{ with } A\in\R^{4\times 4}, 
\mbox{ and } A = A^T.
\ee   
We do not require regularity of $A$, and will later on in Section~\ref{pedal_conics:sec}
discuss singular dual quadrics. In case that $\det A\not=0$, and considering
$X=(x_0,x_1,x_2,x_3)$ as homogeneous point coordinates, the quadric $\ol F$,
whose tangent planes satisfy~\eqref{dualquadric:eq}, is the zero-set 
of the quadratic polynomial
\bes
\ol F(X) &=& X^T \cdot A^{-1} \cdot X.
\ees   

\subsubsection{Pedal surfaces of regular quadrics}

Consider a regular dual quadric $\ol F^\star$, represented as zero-set of a quadratic 
polynomial~\eqref{dualquadric:eq}. Let $O\notin F$ be the reference point of the 
foot-point map. There are two main cases to be distinguished. 
A quadric $\ol F$ is called a paraboloid, if the ideal plane $\omega: x_0=0$
is tangent to $\ol F$. Otherwise, $\ol F$ is
either a hyperboloid or an ellipsoid. Let us start with the latter case.
We may assume that the dual quadric $\ol F^\star$ is represented by
\bes
\ol F^\star: 
a_0u_0^2 + a_1u_1^2 + a_2u_2^2 + a_3u_3^2 + 
u_0(b_1u_1 + b_2u_2 + b_3u_3)          = 0. 
\ees
This choice is justified because we may choose $\omega$, and the symmetry planes of the tangential cone with vertex at $O$ as coordinate planes.
Applying the foot-point map $\alpha$ and with the abbreviation $\xv=(x_1,x_2,x_3)$, 
its pedal surface $\ol G = \alpha(\ol F^\star)$ is represented by 
\be\label{cyclide1:eq}
\ol G: x_0^2(a_1x_1^2 + a_2x_2^2 + a_3x_3^2) 
- x_0(\xv^2)(b_1x_1+b_2x_2+b_3x_3) + a_0 (\xv^2)^2= 0.
\ee
According to Proposition~\ref{footpoint_composition}, the foot-point map $\alpha$ is 
the product $\alpha=\sigma\circ\pi$ of the polarity $\pi$ and the inversion $\sigma$. 
Since $\pi(\ol F^\star)$ is a quadric, $\ol G = \sigma(\pi(\ol F^\star))$ is 
the image of a quadric with respect to the inversion $\sigma$.
Thus $\ol G$ is a special instance of a Darboux cyclide, an algebraic surface 
typically of degree four, whose intersection with $\omega$ is the conic 
$j: x_0=0, x_1^2+x_2^2+x_3^2=0$, with multiplicity two. 
Note that not all Darboux cyclides are images of quadrics with respect to inversion.
Since the highest power of $x_0$ in~\eqref{cyclide1:eq} is two, the 
reference point $O$ is a double point of $\ol G$. We mention that parabolic 
Darboux cyclides are algebraic surfaces of order three, whose intersection with
$\omega$ contains, besides $j$, a real line, 
compare equation~\eqref{parabolic_cyclide:eq}. 

\begin{remark}
A surface is called Darboux cyclide, if it is the zero-set of a quadratic 
polynomial $Q(\yv): \yv^T\cdot B \cdot \yv$ in {\sf pentaspherical coordinates}
$\yv=(y_0, y_1,\ldots, y_4)$, with $y_0^2 = y_1^2+y_2^2+y_3^3+y_4^2$,
which is just $S^3: \yv^T\cdot\diag(-1,1,1,1,1)\cdot\yv=0$. A cyclide
is obtained by applying a stereographic projection $S^3\to\R^3$ to 
the intersection $Q\cap S^3$. This projection is realized by the equations
\bes
x_0=y_0-y_4, x_1=y_1, x_2=y_2, x_3 = y_3.
\ees
Substituting these relations into~\eqref{cyclide1:eq}, and taking into account
that $(y_0-y_4)(y_0+y_4)=y_1^2+y_2^2+y_3^2$, gives the quadratic equation
\bes
a_0(y_0+y_4)^2 - (y_0+y_4)(b_1y_1+b_2y_2+b_3y_3) + (a_1y_1^2+a_2y_2^2+a_3y_3^2)=0,
\ees
representing $Q$. Thus,~\eqref{cyclide1:eq} is a Darboux cyclide.
\end{remark}

We turn to the case where $\ol F^\star$ is a paraboloid. Since $\omega$ is a tangent
plane of $F$, we may choose a coordinate system such that $F$ is parameterized by 
$\fv(u,v) = (u,v,au^2+bv^2+c)$, with $abc\not=0$. If $a$ or $b$ are zero, 
$F$ is a parabolic cylinder, and $\alpha(F)$ is a curve. If $c=0$, $F$ contains $O$. 
Otherwise, $\ol F$ is the zero-set of the 
polynomial $\ol F(X)=ax_1^2+bx_2^2+cx_0^2-x_0x_3$. Its dual equation reads
\bes
\ol F^\star: -4ab u_0u_3 + bu_1^2 + au_2^2 - 4abc u_3^2 = 0.
\ees
The polarity $\pi$ from~\eqref{polarity:eq} maps $\ol F^\star$ to a quadric 
$\pi(\ol F^\star): 4ab x_0x_3 + bx_1^2 + ax_2^2 - 4abc x_3^2 = 0$. Since
$\omega$ is tangent to $\ol F$, the origin $O=(1,0,0,0)\R$ is contained in
$\pi(F^\star)$. Consequently, after canceling out the factor $x_0$, 
the pedal surface
\be\label{parabolic_cyclide:eq}
\ol G = \alpha(\ol F^\star): 
4ab x_3( x_1^2 + x_2^2 + x_3^2) + x_0( b x_1^2 + a x_2^2 - 4abc x_3^2) = 0
\ee
is a parabolic Darboux cyclide, an algebraic surface of degree three.
In case that $F$ is a paraboloid of revolution, thus $a=b$, and additionally
$O$ coincides with the focal point of $F$, thus $c=-\frac{1}{4a}$, 
this surface is reducible and reads $\ol G: (x_1^2+x_2^2+x_3^2)(4ax_3+x_0)=0$.
The first factor $x_1^2+x_2^2+x_3^2=0$ defines the isotropic cone with vertex at $O$,
the second factor $4ax_3+x_0=0$ is the tangent plane at the vertex $(1,0,0,-1/(4a))\R$ of $F$. 

\begin{proposition}
The pedal surface $\ol G=\alpha(\ol F^\star)$ of a regular ellipsoid or hyperboloid 
$\ol F$ is a Darboux cyclide of degree four, with the reference point $O$ 
as double point. The family of offset surfaces $\ol F_d$ of the quadric $\ol F$ 
is mapped by $\alpha$ to the conchoid surfaces 
$\ol G_d$ of the cyclide $\ol G$. The pedal surface $\ol G$ of a paraboloid $\ol F$ 
is a Darboux cyclide of degree three. Only in case that the reference point $O$ 
coincides with the focal point of a paraboloid of revolution $F$, 
its pedal surface is the tangent plane at $F$'s vertex.     
\end{proposition}

In addition to the geometric properties of pedal surfaces of quadrics, 
we want to discuss quadrics in context with rational offset surfaces. 
As proved in~\cite{pet1}, the offset surfaces
$F_d$ of regular quadrics $F$ admit rational parameterizations. 
The construction is not trivial, and thus we provide an outline. 
Any regular quadric $\ol F\subset \P^3$ is the envelope of a 
one-parameter family of cones of revolution $C(u)$, 
with vertices at a focal conic of $\ol F$.
It is possible to parameterize these cones of revolution $C(u)$ in a way 
that the normal vectors of their tangent planes have rational norm. 
Since offsets of cones of revolution are again cones of revolution, 
the rationality of the norm of $F$'s  normal vectors holds for 
its offsets $F_d$, too. 

To construct the pedal surface $\ol G_d=\alpha(\ol F_d^\star)$ of the offsets 
$\ol F_d^\star$ of a quadric $\ol F^\star$, one computes the family of 
pedal curves $c=\alpha(C)$ of the family of cones of revolution $C$. 
Such a pedal curve $c$ is a rational spherical curve of degree four,
and can be constructed in the following way. Let $C^\perp$ be a cone of revolution
with vertex $O$, consisting of lines perpendicular to $C$'s tangent planes.
Let $S$ be the sphere with diameter $OV$, where $V$ is the vertex of $C$.
Then the pedal curve $c=\alpha(C)$ equals the intersection $S\cap C^\perp$,
and thus $O$ is a double point of $c$.  
Finally, the pedal surface $\ol G_d$ contains a rational one-parameter family of 
rational spherical quartic curves $c$, according to the cones of revolution $C$
enveloping $\ol F_d^\star$. 


\paragraph{Rational offsets of paraboloids and their pedal surfaces}    

Consider a paraboloid $F$, and its offset surfaces $F_d$, we intend to derive 
explicit rational parameterizations of the offsets and their pedal surfaces $G_d$,
with respect to the reference point $O=(0,0,0)$. Let 
\bes
\fv(u,v)=(u,v,\frac{1}{2}au^2+\frac{1}{2}bv^2+c), \mbox{ with } abc\not=0,
\ees
be a parameterization of $F$. The corresponding normal vector reads
$\nv(u,v)$ = $(-au,-bv,1)$. In order to determine a re-parameterization,
so that the resulting normal vector has rational norm, we consider
the system of equations
\be\label{paraboloid_normalvector:eq}
\nv(u,v) = \lambda(s,t)(\cos s\cos t, \sin s \cos t, \sin t)=\lambda(s,t)\mv(s,t).
\ee
The right hand side of this equation represents the unit vector $\mv(s,t)$, scaled by
$\lambda(s,t)$. It is more compact to 
represent $\mv(s,t)$ in terms of trigonometric functions, but we either may use any 
rational parameterization $\mv(s,t)$ of $S^2$, or use the Weierstrass substitutions
$\cos x = (1-y^2)/(1+y^2)$, and $\sin x = 2y/(1+y^2)$, where $y=\tan (x/2)$, to convert
the trigonometric representation~\eqref{paraboloid_normalvector:eq} to a rational one.
The system~\eqref{paraboloid_normalvector:eq} has the obvious solution
\be\label{repar_paraboloid:eq}
\lambda(s,t) = 1/\sin t, u = \frac{-\cos s\cos t}{a\sin t}, v=\frac{-\sin s \cos t}{b\sin t}.
\ee
Performing this reparameterization, the tangent planes 
$E$ and $E_d$ of $F$ and its offset surface $F_d$ are represented by
\bes
E_d:     \xv\cdot\mv(s,t) &=&  e(s,t) + d, \nonumber\\
\mbox{ with } e(s,t) &=& 
       -\frac{b\cos^2 s \cos^2 t + a\sin^2 s \cos^2 t - 2abc\sin^2 t}{2a b \sin t}.
\ees
The foot-point map $\alpha$~\eqref{footpoint_map:eq} transforms $E_d$ to 
a parameterization $\gv_d(s,t)$ of the pedal surfaces $G_d$. Since $\|\mv\|=1$, 
we obtain the polar representation 
\bes 
\gv(s,t) = (e(s,t)+d)\mv(s,t)
\ees  
of the conchoid surfaces $G_d$ of $G$. Thus the pedal surfaces $G_d$ of the 
offset surfaces $F_d$ of paraboloids $F$ are rational conchoid surfaces
in the sense of Definition~\ref{ratconch:def}.

\subsubsection{Pedal surfaces of conics as singular dual quadrics}
\label{pedal_conics:sec}

When speaking about pedal surfaces of quadrics and their offsets, it is worth studying
dual singular quadrics, too. Consider a dual singular quadric
\be\label{dualsingquadric:eq}
\ol F^\star(U) &=& U^T \cdot A \cdot U, \mbox{ with } A\in\R^{4\times 4},
A = A^T \mbox{ and } \rk A = 3,
\ee 
which is formed by the planes passing through the tangent lines of a conic $c=\ol F^\star$. 
The null-space of $A$ is the carrier plane $\gamma$ of $c$. To study the pedal surface $\alpha(c)$,
we use the composition~\eqref{footpoint_composition} of $\alpha=\sigma\circ\pi$, with the 
polarity $\pi:{\P^3}^\star\to \P^3$ and the inversion $\sigma:\P^3\to\P^3$, both with respect
to the unit sphere $S^2$. 

The polar image $Q=\pi(c)$ is a quadratic cone, with vertex $V=\pi(\gamma)$. The generating lines
in $Q$ correspond to the pencils of planes in $c$. The cone $Q$ becomes a cylinder in case that
$O\in \gamma$. In case that $O\in c$, a pencil of planes is mapped to points of an ideal line,
and $Q$ is a parabolic cylinder. Worth to be mentioned is also the case that $Q=\pi(c)$ is a 
rotational cone or cylinder. This happens exactly if $O$ is contained in the focal conic $d$
of $c$. The carrier plane of $d$ is a symmetry plane of $\gamma$, and the vertices of $d$
coincide with the focal points of $c$, and vice versa.    

The inversion $\sigma$ maps a quadratic cone $Q$ to a cyclide $G=\sigma(Q)$, 
typically of degree four. 
The family of tangent planes $\tau(u)$ of $Q$ is mapped to a family of spheres 
$S(u)=\sigma(\tau(u))$,
so that $G$ is a canal surface. In case that $Q$ is a cylinder or cone of revolution, 
it contains also
a family of inscribed spheres. Consequently, $G=\sigma(Q)$ is the envelope of 
two different families of spheres, thus a Dupin cyclide.

\begin{proposition}
The pedal surface $G$ of a conic $c$ is a Darboux cyclide, being also a canal surface.
The pedal surfaces of the offset surfaces of $c$, being pipe surfaces with center curve $c$,
are conchoid surfaces $G_d$ of the Darboux cyclide $G$. In case that the reference point $O$
of the foot-point map is located at the focal conic $d$ of $c$, the pedal surface $G$ is a 
Dupin cyclide, being the envelope of a family of spheres in a twofold way. 
\end{proposition}

Typically the pedal surfaces of conics are canal surfaces of degree four. 
Degree reductions appear at first if $c=\ol F^\star$ is a parabola. 
Since $\omega\in\ol F^\star$, $O\in\pi(c)$, and in general
$\alpha(c)$ is a cubic cyclide. This is illustrated in the following.

\paragraph{Example} We consider the parabola $F$ with parameterization
$\fv(u)=(u,0,\frac{a}{2}u^2+c)$, with $ac\not=0$, whose tangent planes are 
$E(u,v)=\R(-\frac{a}{2}u^2+c, au, v, -1)$. 
In terms of homogeneous plane coordinates $U=(u_0,\ldots, u_3)$, the parabola
is the zero set of the polynomial 
\bes
\ol F^\star(U) = u_1^2 -2au_0u_3 - 2ac u_3^2. 
\ees
The pedal surface $\ol G = \alpha(F^\star)$ is a parabolic Darboux cyclide, 
and its defining polynomial reads
\be\label{par_cyclide:eq}
\ol G(X) = x_0(x_1^2 - 2acx_3^2) + 2ax_3(x_1^2+x_2^2+x_3^2).
\ee
We note that applying the polarity $\pi$ to $\ol F^\star$ yields 
the cylinder $\pi(\ol F^\star)(X) = x_1^2-2acx_3^2+2ax_0x_3$. The inversion $\sigma$
maps this surface to the pedal surface~\eqref{par_cyclide:eq}. In case that 
$\pi(\ol F^\star)$ is a cylinder of revolution, it is the envelope of
two different one-parameter families of spheres, including the family of tangent planes.
Thus, for $2ac=-1$ expressing that the origin is the focal point $F$, 
the pedal surface $\ol G$ is a parabolic Dupin cyclide. 
 
Since the family of offset surfaces $F_d^\star$ is mapped by $\alpha$ to a family
of conchoid surfaces $G_d$, we have a look at rational offset parameterizations
of $F^\star$. Therefore we use the reparameterization~\eqref{repar_paraboloid:eq},
which maps $E(u,v)$ to
\bes
E_d(s,t) = \R\left(
\frac{\cos^2 s\cos^2 t - 2ac\sin^2 t}{2a\sin t} + d, 
\cos s \cos t, \sin s \cos t, \sin t\right).
\ees
In fact we use trigonometric instead of rational functions because of readability.  
The offsets of the parabola $\ol F^\star$ are pipe-surfaces $\ol F_d^\star$ 
of class four, whose center curve is $F$. Their dual equation reads
\bes
\ol F_d^\star(U) &=& -4a^2 d^2 u_3^2 (u_1^2+u_2^2 + u_3^2)
+(u_1^2 -2au_0u_3 - 2ac u_3^2)^2, \nonumber\\
&=& -4 a^2 d^2 u_3^2 (u_1^2+u_2^2 + u_3^2) + \ol F^\star(U)^2.
\ees
The conchoid surfaces $\ol G_d$ of the parabolic cyclide $\ol G$ are the $\alpha$-images
of $\ol F_d^\star$. Consequently, their implicit equation reads
\bes
\ol G_d(X) &=& -4a^2 d^2 x_0^2 x_3^2(x_1^2+x_2^2+x_3^2)
+(x_0 (x_1^2 - 2 ac x_3^2) + 2 a x_3 ( x_1^2 + x_2^2 + x_3^2))^2, \nonumber\\
&=& -4a^2 d^2 x_0^2 x_3^2(x_1^2+x_2^2+x_3^2) + \ol G(X)^2.
\ees 

\begin{figure}[ht]
\subfigure[Offset of Parabola]{\label{parabel_offset:fig}
\begin{overpic}[width=.3\textwidth]{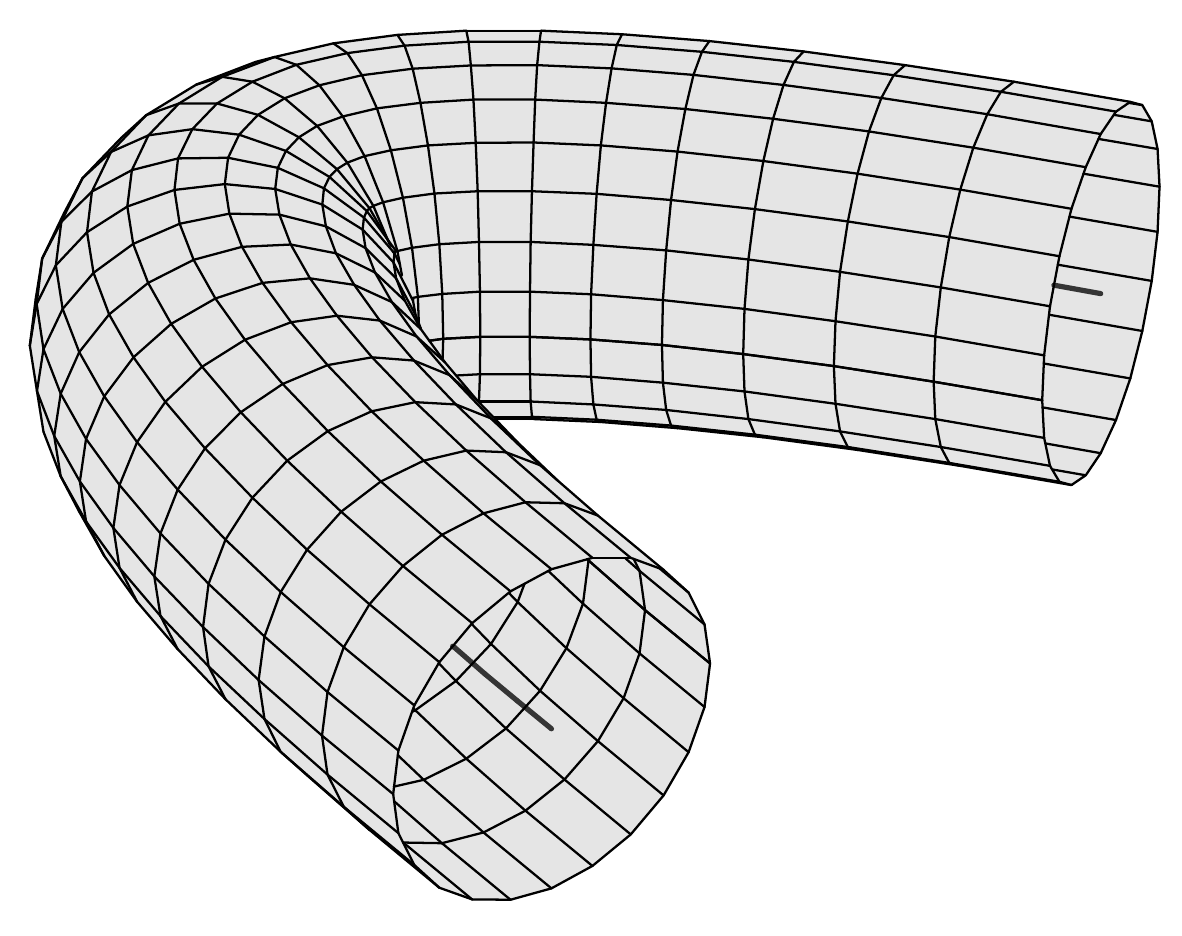}%
\end{overpic}
}
\hfil
\subfigure[Dupin cyclide]{\label{dupin_cyclide:fig}
\begin{overpic}[width=.3\textwidth]{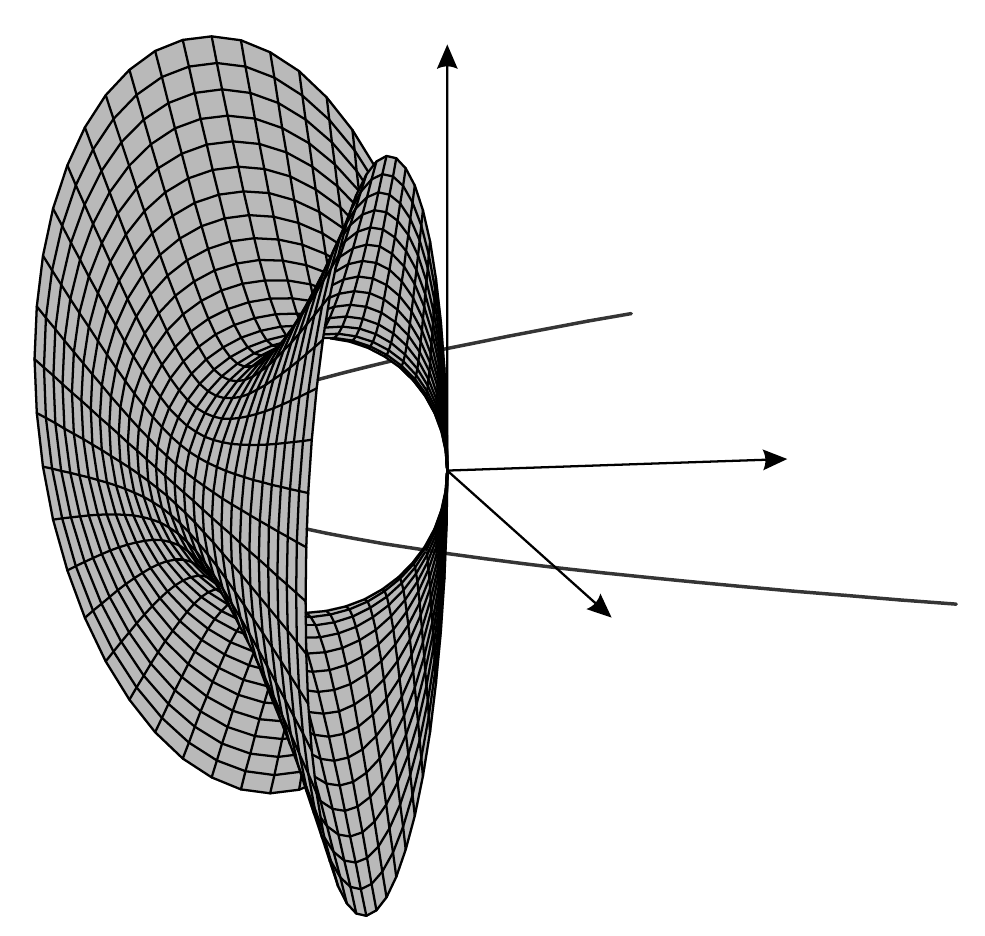}%
\end{overpic}
}
\subfigure[Conchoid of Dupin cyclide]{\label{cyclide_conchoid:fig}
\begin{overpic}[width=.3\textwidth]{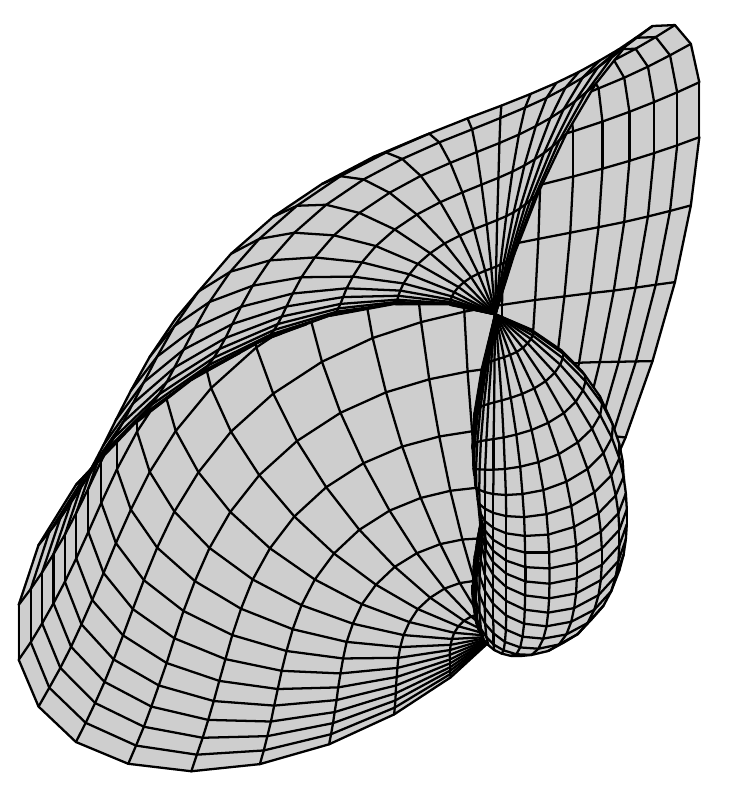}%
\end{overpic}
}
\caption{Parabola $F$, Dupin cyclide $G$ as pedal surface with respect to the parabolas focal point, and conchoid of $G$}
\label{pedalsurf_parabola:fig}
\end{figure}

\paragraph{Quadrics of revolution}
What we have discussed is far from being a complete classification
of the pedal surfaces of quadrics. One important case not being touched so far
is the case where $\pi(\ol F^\star)$ a rotational quadric. This appears if
the tangential cone of the reference point $O$ is a cone of revolution.
In this case, $\pi(\ol F^\star)$ is the envelope of a one-parameter
family of spheres. 

Let $\pi(\ol F^\star)$ be a rotational quadric, then $\ol G$ is the 
envelope of a one parameter family of spheres, since $\sigma$ maps 
spheres to spheres. Thus $\ol G$ is a canal surface. As we have already seen,
this happens also in case that $\pi(\ol F^\star)$ is a quadratic cone or cylinder.
In the particular case where $\pi(\ol F^\star)$ is a cone or cylinder of revolution, 
it is the envelope of two different families of spheres, where tangent planes are counted
as spheres. This implies that $\ol G$ is a Dupin cyclide.

A further particular case appears when $\pi(\ol F^\star)$ is a sphere. 
This happens if the reference point $O$ is the focal point of a rotational 
quadric $\ol F^\star$. Consequently, $\ol G$ is a sphere.

\subsection{Inverse pedal surfaces of quadrics}

While the last section has studied quadrics as dual objects, we consider quadrics $G$ and their
conchoid surfaces $G_d$ as set of points. Let $X=(x_0,x_1,x_2,x_3)$ be homogeneous point 
coordinates in $\P^3$, a quadric $\ol G$ is represented as zero-set
of a quadratic polynomial
\bes
\ol G(X) &=& X^T \cdot A \cdot X, \mbox{ with } A\in\R^{4\times 4}, 
\mbox{ and } A = A^T.
\ees
If $\det A\not=0$, the quadric is regular, and in case that $\rk A = 3$, the quadric
is a quadratic cone. More degenerated cases are a pair of planes ($\rk A=2$) and 
a two-fold plane ($\rk A=1$).
The inverse pedal surface $F^\star = \alpha^\star(G)$ of a plane $G$ is a paraboloid
of revolution, with focal point $O$. 
For the remainder of this section we assume that $\rk A = 3$ or $4$. 

Consider a regular quadric $\ol G$ and the reference point $O$ as origin.
Assuming that the axes of $\ol G$ are parallel to the coordinate axes, we may assume 
that $\ol G$ is the zero-set of the quadratic polynomial
\bes
\ol G(X) = a_0x_0^2 + a_1x_1^2 + a_2x_2^2 + a_3x_3^2 
+ x_0(b_1x_1 + b_2x_2 + b_3x_3).
\ees 
Applying $\alpha^\star$ yields $\ol F^\star$, which is in fact just the dual object of~\eqref{cyclide1:eq}, and its defining polynomial in homogeneous plane coordinates $U=(u_0,u_1,u_2,u_3) = (u_0, \uv)$ reads
\be\label{dualcyclide1:eq}
\ol F^\star(U)= u_0^2(a_1u_1^2 + a_2u_2^2 + a_3u_3^2)
-u_0\uv^2(b_1u_1 + b_2u_2 + b_3u_3) + a_0 (\uv^2)^2.
\ee 
Dual to the fact that the origin $O=(1,0,0,0)\R$ is a double point of $\ol G$ 
in~\eqref{cyclide1:eq}, the ideal plane $\omega$ is 
a double tangent plane of $\ol F^\star$.

In~\cite{quadric_conchoid} it is shown that quadrics $G\subset\R^3$ and 
their conchoid surfaces $G_d$ allow rational parameterizations 
$\gv_d$ with rational norm $\|\gv_d\|$, independently on the chosen reference point.    
According to Theorem~\ref{corr-theo:theo}, rational parameterizations $\fv_d$ of $\ol F_d^\star$
are obtained by $\fv_d = \alpha^\star(\gv_d)$. Using the composition $\alpha^\star=\pi^\star\circ\sigma$,
the surface $\ol F^\star$ is the dual object of a cyclide $\sigma(\ol G)$, and the offset
surface $\ol F_d^\star$ is the dual object of a conchoid surface $\sigma(\ol G_d)$ of a cyclide.

\paragraph{Inverse Pedal surfaces of quadratic cones and cylinders}

Consider a quadratic cone or cylinder $G$, as singular quadric. 
To get insight to the geometric properties of their inverse pedal
surfaces we use the decomposition $\alpha^\star=\pi^\star\circ\sigma$.
The surface $\sigma(G)$ is a cyclide, and since $G$ is the envelope
of a one-parameter family of tangent planes $\tau$, $\sigma(G)$ is
a canal surface, enveloped by the spheres $\sigma(\tau)$. These spheres
$\sigma(\tau)$ pass through $O$. In case $O\in\tau$, $\sigma(\tau)=\tau$ is a plane
as well. 

Quadratic cones and cylinders are also ruled surfaces, and thus we can use the 
results from Section~\ref{invpedal_rs:sec}. There we noted that 
$\alpha^\star(G) = F^\star$ is the envelope of parabolic cylinders 
$P=\alpha^\star(l)$, being the $\alpha^\star$-images of the generating lines $l\subset G$. 

The quadrics $\alpha^\star(\tau)$ enveloping $F$, are typically paraboloids,
according to the fact that the spheres $\sigma(\tau)$ contain $O$, and that the polarity
$\pi^\star$ maps a sphere $\sigma(\tau)\ni O$ to a paraboloid. Since $l\subset \tau$,
the parabolic cylinders $\alpha^\star(l)$ touches the paraboloids $\alpha^\star(\tau)$
in points of conics $d$. 
Consider the cones of revolution $D$ touching the canal surface $\sigma(G)$ 
along its characteristic circles $\sigma(l)$. The previously mentioned conics $d$
are just the images of these cones $D$ with respect to $\pi$, thus $d=\pi(D)$. 

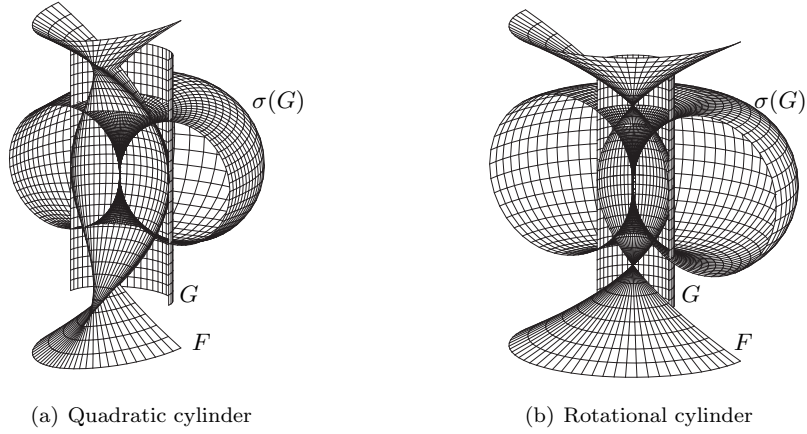
\begin{figure}[ht]
\centering
\subfigure[Quadratic cylinder]{\label{pedal_cylinder1:fig}
\begin{overpic}[width=.3\textwidth]{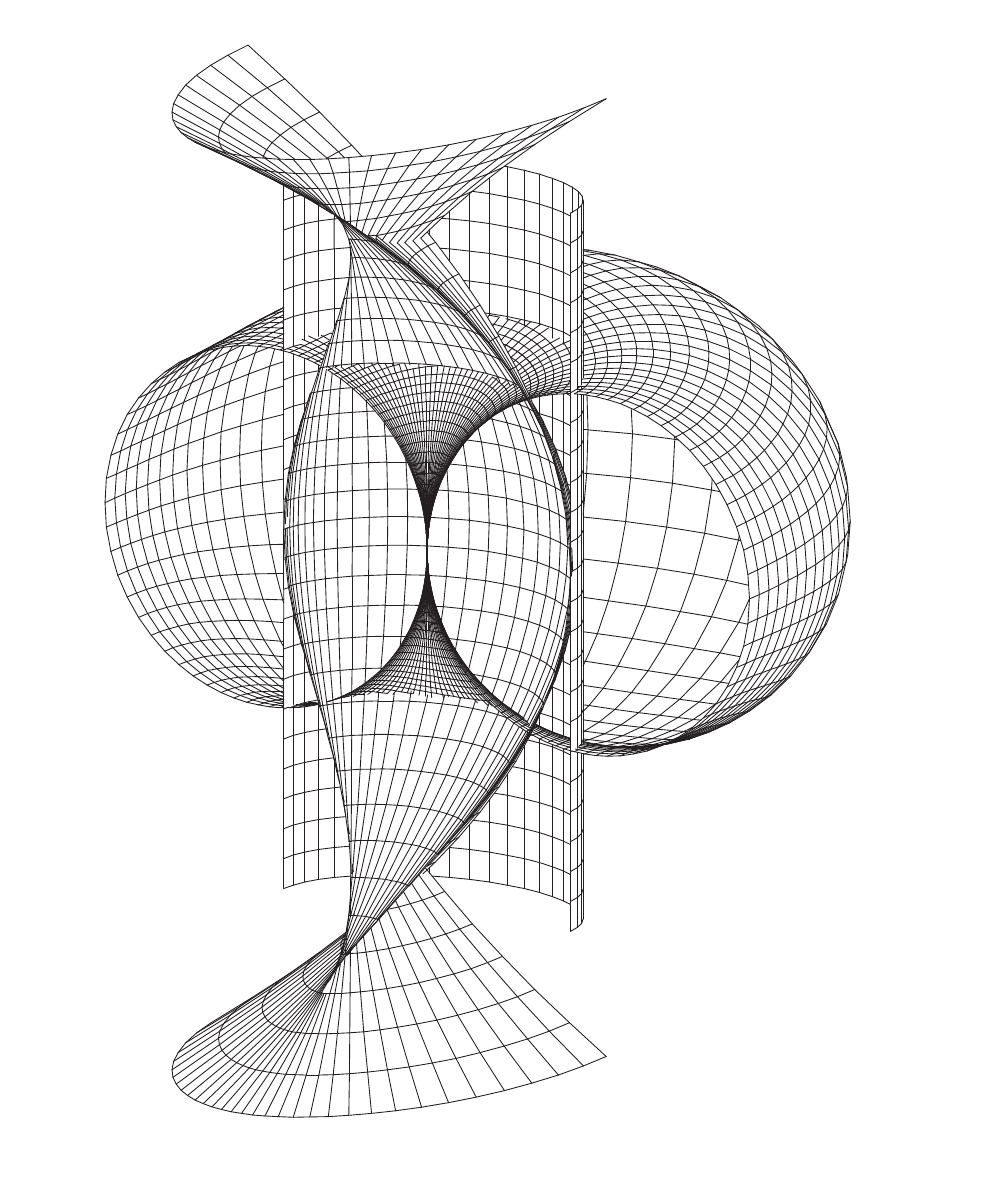}%
{\small
\put(51,22){$G$} \put(69,70){$\sigma(G)$} \put(54,10){$F$} 
}
\end{overpic}
}
\hfil
\subfigure[Rotational cylinder]{\label{pedal_cylinder2:fig}
\begin{overpic}[width=.3\textwidth]{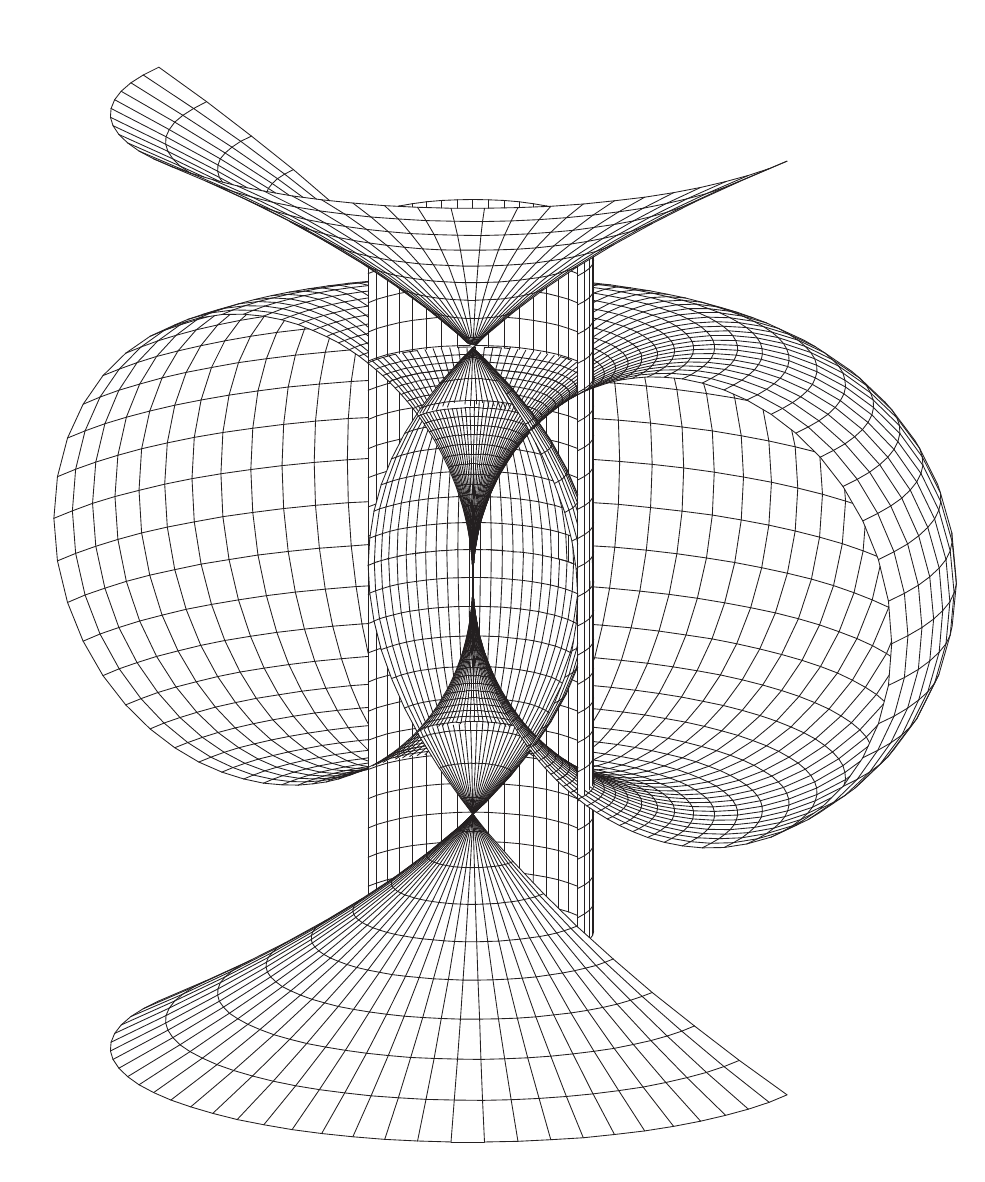}%
{\small
\put(52,22){$G$} \put(70,70){$\sigma(G)$} \put(65,10){$F$} 
}
\end{overpic}
}
\caption{Inverse pedal surfaces of quadratic cylinders}
\label{invpedalcylinder:fig}
\end{figure}

We illustrate this at hand of an example.
Consider the cylinder $G: x^2/a^2 + y^2/b^2-1=0$. In terms of homogeneous coordinates
it reads $\ol G: x_1^2/a^2 + x_2^2/b^2-x_0^2=0$.
The inverse image $\sigma(\ol G)$ is a cyclide. Applying $\pi: \P^3\to {\P^3}^\star$, 
one obtains $\alpha^\star(\ol G) = \ol F^\star$. These surfaces read
\bes
\sigma(\ol G) &:& \frac{x_0^2x_1^2}{a^2} + \frac{x_0^2x_2^2}{b^2}=(x_1^2+x_2^2+x_3^2)^2,\\
\ol F^\star &:&   \frac{u_0^2u_1^2}{a^2} + \frac{u_0^2u_2^2}{b^2}=(u_1^2+u_2^2+u_3^2)^2.
\ees 
The defining polynomial of $\ol F$ is rather lengthy, and of degree eight. Thus we provide 
a parametric representation. Starting with $\gv(u,v)=(a \cos u, b \sin u, v)$ for $G$, 
the tangent planes of $F$ are $\alpha^\star(\gv)$ $=$ 
$E=\R(-a^2 \cos^2 u - b^2\sin^2 u - v^2, a\cos u, b\sin u, v)$. 
An affine parameterization of $F$, see Figure~\ref{pedal_cylinder1:fig}, is finally
\bes
\fv(u,v) = 
(-\frac{\cos u}{a}((a^2-b^2) \cos^2 u + b^2 - 2a^2 + v^2), 
       -\frac{\sin(u)}{b}((a^2-b^2)\cos^2 u - b^2 + v^2), 2v). 
\ees 
Even the case $a=b$ is interesting, where $\sigma(\ol G)$ is a torus, whose 
meridian circles touch the $z$-axis. The inverse pedal surface is a rotational surface 
with a parabola $(-b^2+v^2,0,2v)$ as meridian curve, see Figure~\ref{pedal_cylinder2:fig}.

To represent the whole family of conchoid surfaces $G_d$ of the quadratic cylinder $G$, 
one preferably uses a rational polar representation, see~\cite{pet10}. Applying 
$\alpha^\star$ gives a dual representation of the family $\ol F_d$ of rational
offset surfaces. Let $\dv=(d_1,d_2,d_3)$ $=$ $(a\cos u, b\sin u, 0)$ be the foot-point curve
with respect to $O$, a polar representation of $G_d$ reads
\bes
\gv_d(u,t) = \frac{1+\dv^2t^2+2dt}{2t(1+\dv^2t^2)}(2t d_1, 2t d_2, 1-\dv^2t^2).
\ees 
The tangent planes of $F_d$ are $E_d:\xv\cdot\gv_d = \gv_d^2$. Solving the linear
system $E_d\cap (E_d)_t \cap (E_d)_u$ yields an affine parameterization $\fv_d(u,v)$ 
of $F_d$, whose normal vector $\gv_d$ has rational length.

\begin{figure}[ht]
\centering
\subfigure[Ellipse]{\label{circle2offset1_fig}
\begin{overpic}[width=.3\textwidth]{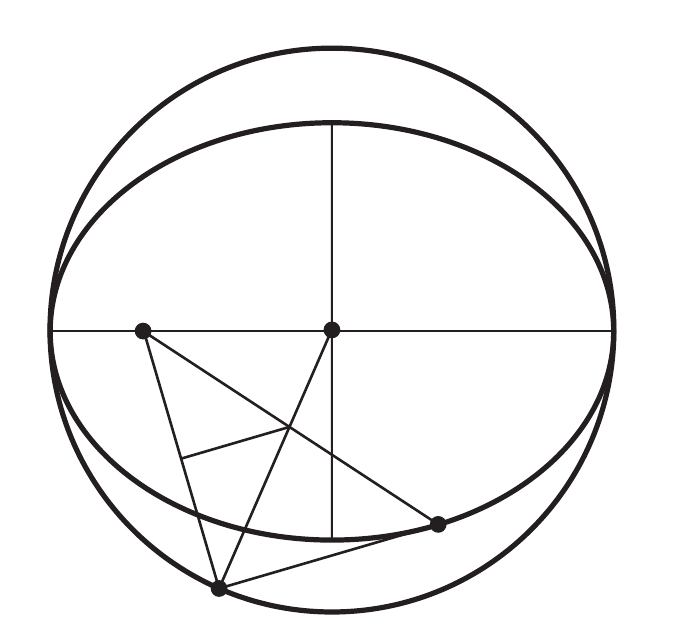}%
{\small
\put(51,49){$M$}  \put(22,49){$O$} \put(40,35){$Y$}
\put(29,0){$X$} \put(60,10){$X'$} \put(30,30){$b$}
\put(60,85){$C$} \put(60,64){$C'$}
}
\end{overpic}
}
\hfil
\subfigure[Hyperbola]{\label{circle2offset2_fig}
\begin{overpic}[width=.3\textwidth]{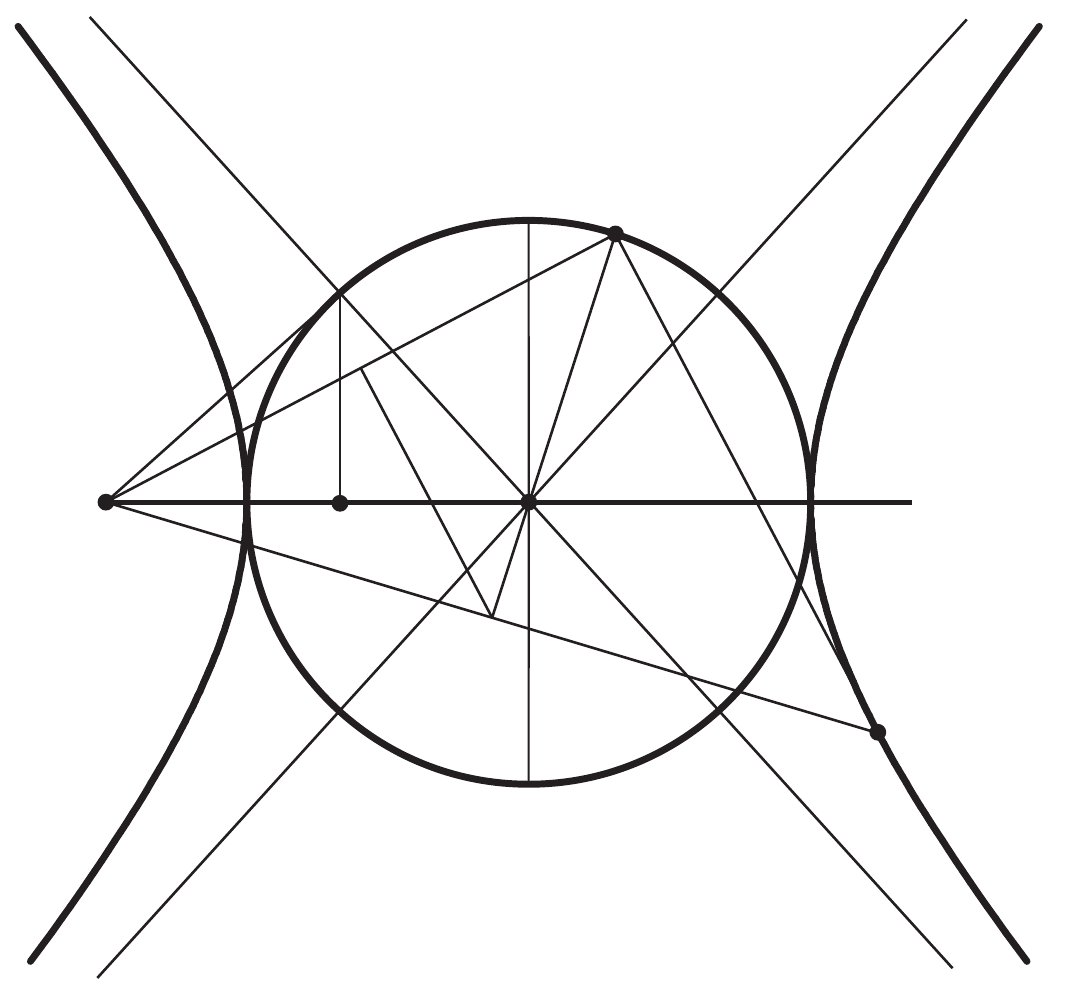}%
{\small
\put(54,45){$M$}  \put(5,38){$O$} \put(32,38){$O'$} \put(42,28){$Y$}
\put(56,72){$X$} \put(85,22){$X'$} \put(39,46){$b$}
\put(40,72){$C$} \put(87,72){$C'$}
}
\end{overpic}
}
\caption{Inverse pedal curves $C'$ of a circle $C$ with respect to point $O$}
\label{invpedalcircle:fig}
\end{figure}
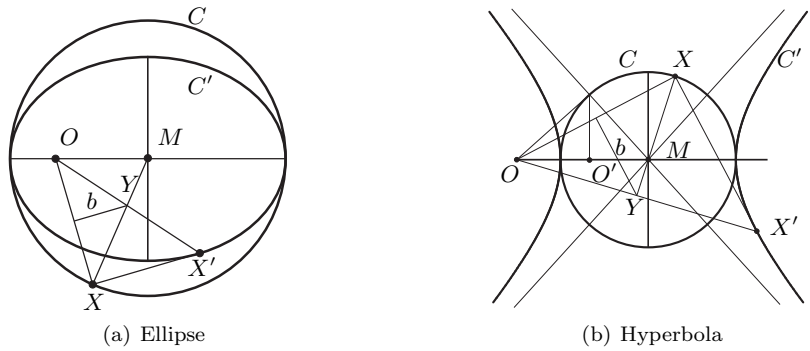

\paragraph{Inverse pedal surface of a sphere}
    
As a particular example we discuss the inverse pedal surface $F^\star$ of a sphere 
$G:(x-m^2)+y^2+z^2=r^2$, of radius $r$ and centered at $(m,0,0)$, with 
respect to the reference point $O=(0,0,0)$. Since $\sigma(G)$ is a sphere,
$F^\star = \pi^\star(\sigma(G))$ is a quadric.
The defining polynomials of $\ol G$ and $\ol F^\star$ in  
homogeneous coordinates are
\bes
\ol G(X) &=& x_1^2+x_2^2+x_3^2 - 2x_0x_1m + x_0^2(m^2-r^2),\\
\ol F^\star(U) &=& u_0^2 + 2u_0u_1m + (u_1^2+u_2^2+u_3^2)(m^2-r^2). 
\ees 
Transforming plane to point coordinates one obtains the known result that
$\ol F$ is a quadric of revolution, with $O$ as focal point,
\bes
\ol F:  \frac{r^2(y^2+z^2)}{a^2} - \frac{x^2}{a} + \frac{2mx}{a} = 1, 
\mbox{ with } a=m^2-r^2.
\ees
The quadric $\ol F$ is an ellipsoid for $O$ being inside of $\ol G$, and a
hyperboloid of two sheets for $O$ being outside of $\ol G$. 
Figure~\ref{invpedalcircle:fig} illustrates the 2d-case, 
the inverse pedal curves of a circle with respect to a reference point $O$.
For $O\in\ol G$, the inverse pedal surface $\ol F$ degenerates to a single point, 
according to Thales theorem.
Considering the dual surface in that particular case, 
$\ol F^\star$ is a bundle of planes, passing through that 
single point. 

According to Theorem~\ref{corr-theo:theo}, a particular parameterization $\gv(u,v)$
of $G$ with rational norm $\|\gv(u,v)\|$, see \cite{pet12} is mapped by $\alpha^\star$ 
to a parameterization $\fv(u,v)$ of $F$, whose normal vector has rational length. 
In this way, the conchoid surfaces of spheres correspond to the offsets of ellipsoids 
and hyperboloids of revolution.   

\begin{corollary}
The conchoid surfaces $G_d$ of a sphere $G$ are in birational correspondence to the 
offset surfaces $F_d$ of ellipsoids or hyperboloids of revolution $F$, 
depending on whether the reference point $O$ is inside or outside of $G$.
\end{corollary}

\section{Conclusion}

The present article studies relations between families of rational offset surfaces
and rational conchoid surfaces. The foot-point map $\alpha$ transforms a family of offset
surfaces to a family of conchoid surfaces, where the reference points of the foot-point
map and the conchoid construction coincide. Since $\alpha$ is a birational map, 
birational invariants are transformed from one family to the other. The same properties
hold for the inverse foot-point map. 
The relations between offset surfaces and conchoid surfaces are demonstrated
at hand of pedal surfaces and inverse pedal surfaces of ruled surfaces and quadrics. 

There is a close relation to bisector surfaces. Considering a surface $G\subset \R^3$,
and a fixed reference point $O=(0,0,0)$. The bisector surface $B(G,O)$ of $G$ and $O$ 
is the envelope of symmetry planes $S$ of $O$ and a moving point $\gv\in G$. 
Scaling the inverse pedal surface $F^\star = \alpha^\star(G)$ by the factor $1/2$
gives $B(G,O)$. More details on this construction may be found in~\cite{pet_bisector}.

\section*{Acknowledgments}
This work has been partially funded
by the Spanish 'Ministerio de Economia y Competitividad'  under the project
MTM2011-25816-C02-01.
Last two authors belongs to  the Research Group ASYNACS (Ref. CCEE2011/R34).



\begin{thebibliography}{}


\bibitem{alberto}
A.~Albano,  M.~Roggero: Conchoidal transform of two plane curves,
Applicable Algebra in Engineering, Communication and Computing,
Vol.21, No.4, 2010, pp. 309--328.

\bibitem{arrondo}
Arrondo, E., Sendra, J., Sendra, J.R., 1997.
Parametric Generalized Offsets to
Hypersurfaces. Journal of Symbolic Computation {\bf 23}, 267–-285.




\bibitem{dietz}
Dietz, R., Hoschek, J., and J\"uttler, B., 1993. 
An algebraic approach to curves and surfaces on the sphere
and other quadrics, Comp. Aided Geom. Design {\bf 10}, 211–229.

%
%

\bibitem{F1} Farouki, R.T., 2002. Pythagorean Hodograph curves, in: Handbook of Computer Aided
Geometric Design, Farin, G., Hoschek, J., Kim, M.-S. eds.), Elsevier, 2002.

\bibitem{F2} Farouki, R.T., 2008.
Pythagorean-Hodograph Curves: Algebra and Geometry Inseparable, Springer.

\bibitem{quadric_conchoid}
Gruber, D. and Peternell, M.: Conchoid surfaces of quadrics, submitted
to JSC.

\bibitem{KPet}
R. Krasauskas,  M. Peternell, 2010. 
Rational offset surfaces and their modeling applications, in:
IMA Volume 151: Nonlinear Computational Geometry, (eds.) I.Z. Emiris, F. Sottile, and Th.
Theobald, 2010, pp. 109–135.



\bibitem{HL}
J. Hoschek, D. Lasser (1993). Fundamentals of Computer Aided Geometric
Design. A.K. Peters, Ltd. Natick, MA, USA (1993)

\bibitem{pet1} Peternell, M. and Pottmann, H., 1998.
    A Laguerre geometric approach to rational
    offsets, Comp.\ Aided Geom.\ Design {\bf 15},
    223--249.
%
%
\bibitem{pet_bisector} Peternell, M., 2000.
Geometric Properties of Bisector Surfaces
Graphical Models and Image Processing {\bf 62}, 202--236.

\bibitem{pet10}
Peternell, M., Gruber, D. and Sendra, J., 2011:
Conchoid surfaces of rational ruled surfaces,
Comp.\ Aided Geom.\ Design {\bf 28}, 427--435.%

\bibitem{pet12}
Peternell, M., Gruber, D. and Sendra, J., 2012:
Conchoid surfaces of spheres,
Comp.\ Aided Geom.\ Design {\bf 30}, issue 1, pp. 35-44. %
.%
%
\bibitem{pot98} Pottmann, H., and Peternell, M. 1998.
   Applications of Laguerre geometry in CAGD,
   offsets, Comp.\ Aided Geom.\ Design {\bf 15},
   165--186.
%

\bibitem{NSW} Ngô L.X.C., Sendra J.R., Winkler F. 2012. 
Birational Transformations on Algebraic Ordinary Differential Equations.
With Ngô L.X.C., Winkler F.
RISC report, no. 12-18 (2012)


\bibitem{SSS} San Segundo F., Sendra J.R. 2012. Total Degree Formula for the Generic Offset to a Parametric Surface
International Journal of Algebra and Computation Vol. 22, No. 2.

%
\bibitem{schicho1} Schicho, J., 2000.
Proper Parametrization of Real Tubular Surfaces,
J. Symbolic Computation {\bf 30}, 583--593.
%


\bibitem{sendra0} J.R. Sendra, J. Sendra, 1999. Algebraic Analysis of Offsets to Hypersurfaces.
Mathematische Zeitschrift Vol. 234, 697–719.

\bibitem{sendra1} J.R.~Sendra and J.~Sendra, 2008.
An algebraic analysis of conchoids to algebraic curves,
Applicable Algebra in Engineering, Communication and Computing,
Vol. 19, No.5, pp.~285--308.

\bibitem{sendra2} J.~Sendra and J.R.~Sendra, 2010.
Rational parametrization of conchoids
to algebraic curves, Applicable Algebra in Engineering, Communication and Computing,
Vol.21, No.4, pp.~413--428.
%

\bibitem{sendra3} J.R. Sendra, D. Sevilla, 2011.
Radical Parametrizations of Algebraic Curves by Adjoint Curves
Journal of Symbolic Computation 46, 1030-1038.

\bibitem{sendra4} J.R. Sendra and D. Sevilla, 2013. First Steps Towards Radical Parametrization of Algebraic Surfaces.
Computer Aided Geometric Design (In press).


 \bibitem{libro} J.R. Sendra J.R., F. Winkler, S. P\'erez-Diaz (2007).
 {\it Rational Algebraic Curves: A Computer Algebra Approach}. Springer-Verlag Heidelberg, in series Algorithms and Computation  in Mathematics.  Volume 22.

\bibitem{VL} Vr\v{s}ek J., L\'avicka M., 2012. Exploring hypersurfaces with offset-like convolutions. Computer Aided Geometric Design Vol. 29, Issue 9, pp. 676-690.


\end{thebibliography}
\end{document}